\documentclass[draft]{amsart}
\usepackage{amsmath,amssymb,amsthm,verbatim}
\usepackage[mathscr]{eucal}

\allowdisplaybreaks 
\numberwithin{equation}{section}

\newtheorem{theorem}{Theorem}
\newtheorem*{theorem*}{Theorem}

\newtheorem*{conjecture*}{Conjecture}

\newtheorem{lemma}{Lemma}

\newtheorem{corollary}{Corollary}
\newtheorem*{corollary*}{Corollary}

\newtheorem*{HypoA}{Hypothesis $A$}

{Theorem (Bombieri--Vinogradov for $\Lambda*\Lambda$)}

\newcommand{\Yildirim}{Y{\i}ld{\i}r{\i}m}

\newcommand{\cL}{{\mathcal L}}

\newcommand{\cS}{{\mathcal S}}
\newcommand{\bR}{{\mathbb R}}

\newcommand{\ctS}{\tilde{\cS}}

\newcommand{\fS}{\mathfrak S}
\newcommand{\asymptotic}{\sim}
\newcommand{\binomial}{\binom}

\newcommand{\ds}{\displaystyle}
\newcommand{\li}{\text{li}}

\newcommand{\Pt}{\tilde{P}}

\def\sumprime_#1{\setbox0=\hbox{$\scriptstyle{#1}$}
\setbox2=\hbox{$\displaystyle{\sum}$}
\setbox4=\hbox{${}'
\mathsurround=0pt$}
\dimen0=.5\wd0 \advance\dimen0 by-.5\wd2
\ifdim\dimen0>0pt
\ifdim\dimen0>\wd4 \kern\wd4 \else\kern\dimen0\fi\fi
\mathop{{\sum}'}_{\kern-\wd4 #1}}

\def\sumflat_#1{\setbox0=\hbox{$\scriptstyle{#1}$}
\setbox2=\hbox{$\displaystyle{\sum}$}
\setbox4=\hbox{${}\flat
\mathsurround=0pt$}
\dimen0=.5\wd0 \advance\dimen0 by-.5\wd2
\ifdim\dimen0>0pt
\ifdim\dimen0>\wd4 \kern\wd4 \else\kern\dimen0\fi\fi
\mathop{{\sum}^\flat}_{\kern-\wd4 #1}}

\begin{document}

\title[Small gaps between products of two primes]
{S\lowercase{mall gaps between products of two primes}}

\author{D. A. Goldston, S. W. Graham, J. Pintz and C. Y. \Yildirim}

\date{September 20, 2006}

\maketitle

\section{Introduction}
\label{sec:1}

As an approximation to the twin prime conjecture it was proved
in \cite{GPY1} that
\begin{equation}
\liminf_{n \to \infty} \frac{p_{n+1} - p_n}{\log p_n} = 0.
\label{eq:1.1}
\end{equation}
The strongest approximation for the twin prime conjecture in
another direction was proved  in the celebrated work
of Chen \cite{Che2} 
\footnote{Chen's result was announced in 1966 \cite{Che1}.
However, due to the Cultural Revolution,  the complete proof was
not published until 1973.}
(see also \cite[Chapter~11]{HR}), where
he showed that  there are infinitely many primes $p$ such that
$p + 2\in \mathcal P_2$, where
\begin{equation}
\mathcal P_2 := \{n:  \Omega(n) \leq 2\}.
\label{eq:1.2}
\end{equation}
If $\mathcal P$ denotes the set of primes, then Chen's theorem
asserts that at least one of the relations
\begin{equation}
p + 2 = p' \in \mathcal P
\label{eq:1.3}
\end{equation}
or
\begin{equation}
p + 2 = p_1p_2, \quad p_1, p_2 \in \mathcal P
\label{eq:1.4}
\end{equation}
holds for infinitely many primes~$p$.

The phenomenon that we cannot specify which one of the two
equations \eqref{eq:1.3} and \eqref{eq:1.4} has infinitely many
solutions (in reality most probably both, naturally) is the most
significant particular case of the parity problem, a heuristic
principle 
stating that sieve methods cannot differentiate
between integers with an even and an odd number of prime factors.
This principle is based on some extremal examples of Selberg
(see \cite[Ch.~4]{Gr}, \cite[p. 204]{Sel}).
Accordingly, until very recently, problems involving numbers
that are products of two distinct prime factors (which we called
$E_2$-numbers in \cite{GGPY1}) seemed to be as difficult as
problems involving primes, since sieve methods seemed to be not
suitable to attack these problems due to the parity problem.
For example, the analogue of \eqref{eq:1.1},
\begin{equation}
\liminf_{n\to \infty} \frac{q_{n+1} - q_n}{\log q_n / \log\log
q_n} = 0
\label{eq:1.5}
\end{equation}
(where $q_1 < q_2 < \dots$ denotes the sequence of $E_2$-numbers)
was, similar to \eqref{eq:1.1}, not known.

The present authors observed that the method -- a variant of
Selberg's sieve -- which led to the proof of \eqref{eq:1.1} in
\cite{GPY1}, can be used even more successfully for
$E_2$-numbers.
In our preceding work \cite{GGPY1} we gave an alternative proof
of \eqref{eq:1.1}; further we showed that $E_2$-numbers are
infinitely often a bounded distance apart, more precisely,
\begin{equation}
\liminf_{n \to \infty} (q_{n+1} - q_n) \leq 26.
\label{eq:1.6}
\end{equation}

The relation \eqref{eq:1.6} was actually a simple consequence of
a more general result, according to which every admissible (see
the definition below) 8-tuple contains at least two
$E_2$-numbers infinitely often.
The following far reaching generalization of the twin prime
conjecture was formulated qualitatively 100 years ago by L. E.
Dickson \cite{Dic}, and two decades later in a quantitative form
by Hardy and Littlewood \cite{HL}.
In order to formulate the conjecture we define a set
\begin{equation}
\mathcal H = \{h_i\}^k_{i = 1} \qquad
h_i \in \mathbb Z^+ \cup \{0\}
\label{eq:1.7}
\end{equation}
to be admissible if for every prime number $p$ the set $\mathcal
H$ does not cover all residue classes $\text{\rm mod}\, p$.

\medskip
\noindent
{\bf Prime-tuple conjecture.}
{\it Given any admissible set $\mathcal H$, there are infinitely
many integers $n$ such that all numbers of the form $n + h_i$
$(1 \leq i \leq k)$ are primes.
The number of such $n$'s below $N$ is asymptotically equal to
\begin{equation}
 \frac{N}{\log^k N} \mathfrak S(\mathcal H) = 
 \frac{N}{\log^k N}
 \prod_p \left(1 -\frac{\nu_p(\mathcal H)}{p}\right) \left(1 - \frac1p\right)^{-k},
\label{eq:1.8}
\end{equation}
where $\nu_p(\mathcal H)$ denotes the number of residue classes
$\mod p$ covered by $\mathcal H$.}

\medskip
The above conjecture includes (as the case $k = 2)$ the
generalized twin prime conjecture, which states that 
every even number can
be written as the difference of two primes in infinitely many ways.
This was formulated by de Polignac
\cite{Pol} in 1849 in a qualitative way, and in the same work
of Hardy and Littlewood \cite{HL} in a quantitative form.

If we substitute primes by almost primes of the form $P_r$
(integers having at most $r \ge 2$ prime factors) then the qualitative
form of the analogous conjecture is true for $k = 2$, as shown
by Chen's theorem \eqref{eq:1.2}, even for $r = 2$.
This trivially implies that we have infinitely often at least
two $P_2$-numbers in any admissible $k$-tuple for any $k \geq 2$.

We will examine the problem whether for any $\nu$ we can guarantee
that there are infinitely often at least $\nu$ $P_2$-numbers (or
at least $\nu$ $P_r$-numbers with a given fixed $r$, independent
of $\nu$) in any admissible $k$-tuple if $k$ is sufficiently
large, that is, $k \geq C_0(\nu)$.

Such a result seems to be unknown for any fixed value of~$r$.
The strongest result in this direction is due to Heath-Brown
\cite{HB} who showed that  if $\{h_i\}^k_{i = 1}$ is an 
admissible $k$-tuple then there are infinitely many $n$ such  that 
\begin{equation}
\max_{1\leq i \leq k} \omega(n + h_i) < C \log k.
\label{eq:1.9}
\end{equation}
This improved an earlier result of Halberstam and Richert
\cite[Ch. 10]{HR}, where the analogue of \eqref{eq:1.9} was
proved with the max replaced by the average of $\omega(n + h_i)$.

In the case of the primes it was shown in \cite{GPY1} that if
the level $\vartheta$ of distribution of primes (see the
definition \eqref{eq:1.17} below) is any 
fixed number in $(1/2,1]$,  then
\begin{equation}
\liminf_{n \to \infty} (p_{n+1} - p_n) < \infty.
\label{eq:1.10}
\end{equation}

On the other hand, we needed the Elliott--Halberstam conjecture (EH) 
(see \cite{EH})  in its full strength to obtain
\begin{equation}
\liminf_{n \to \infty} \frac{p_{n+2} - p_n}{\log p_n} = 0.
\label{eq:1.11}
\end{equation}
For $p_{n+3} - p_n$ the best result we were able  to prove on EH
in \cite{GPY3} was
\begin{equation}
\liminf_{n \to \infty} \frac{p_{n+3} - p_n}{\log p_n} \leq
e^{-\gamma} \big(\sqrt{3} - \sqrt{2}\big)^2 .
\label{eq:1.12}
\end{equation}

The incredible depth of the assumption EH in
\eqref{eq:1.11}--\eqref{eq:1.12} suggests that it might be very
difficult to prove
\begin{equation}
\liminf_{n \to \infty} (q_{n+\nu} - q_n) < \infty,
\label{eq:1.13}
\end{equation}
already for $\nu = 2$ or $3$.
We will show, however, that our method can be applied very
efficiently to this problem.

In the present work we will show  the existence of at least $\nu$
$E_2$-numbers in any admissible $k$-tuple if $k \geq C_1(\nu)$.
We will also show that $C_1(2) = 3$ is permitted, that is,
every admissible triplet contains at least two $E_2$-numbers
infinitely often.

The mentioned work of Heath-Brown \cite{HB} is based on a method
of Selberg \cite{Sel}.
Selberg considered only the case $k = 2$ and showed that there
are infinitely many pairs $n$, $n + 2$ such that one of them  is a
$P_2$-number, the other a $P_3$-number.

Our method, a modified form of the above mentioned
methods of Selberg and Heath-Brown, also shows that 
$C_0(2) =2$. So we have 
\begin{equation}
n, n + 2 \in \mathcal P_2,
\label{eq:1.14}
\end{equation}
infinitely often, improving Selberg's result but falling short of  \eqref{eq:1.2}.

We will, in fact, prove the above results in the following more
general form, similar to Heath-Brown \cite{HB}.
Let
\begin{equation} 
L_i(x) = a_i x + b_i \quad (1 \leq i \leq k) \quad
a_i, b_i \in \mathbb Z, \quad a_i > 0
\label{eq:1.15}
\end{equation}
be an admissible $k$-tuple of distinct linear forms.
In other words, we
suppose that for every prime $p$ there exists $x_p \in \mathbb Z$ such that
\begin{equation} 
p \nmid \prod^k_{i = 1} (a_i x_p + b_i).
\label{eq:1.16}
\end{equation}

In order to formulate the results we will introduce the level
$\vartheta$ of distribution of primes in arithmetic progressions.
We say that the primes have level of distribution $\vartheta$ if
for any positive $A$ there exists a constant $C=C(A)$ such that 

\begin{equation} \label{eq:1.17}
\sum_{q\leq N^{\vartheta}(\log N)^{-C}} \max_{\substack{a\\ (a,q)= 1}} 
\left(  \sum_{p \equiv a \pmod q)} 1 -\frac{\li (N)}{\varphi(q)} \right)
   \ll_{A} \, \frac{N}{(\log N)^A}.
\end{equation}

The Bombieri--Vinogradov Theorem states that $\vartheta = 1/2$
is admissible.
Elliott and Halberstam \cite{EH} conjectured that
\eqref{eq:1.17} is true for any $\vartheta < 1$.
Friedlander and Granville\cite{FG} proved that \eqref{eq:1.17} is not true with
$\vartheta=1$, but it is possible that it still holds for any fixed $\vartheta < 1$.

In the following we suppose that an analogue of \eqref{eq:1.17} is true for
$E_2$-numbers  with the same value of~$\vartheta$.
This is  true with $\vartheta = 1/2$ unconditionally,
as shown by  Motohashi \cite{Mot}.  Motohashi gives a more general
result; he proves that if two functions satisfy analogues of the Bombieri-Vinogradov 
Theorem, then under certain reasonable conditions, the convolution of the 
two functions also satisfies  an analogue of Bombieri-Vinogradov. 
This may also be proved using a slight variation of the argument of 
Bombieri \cite[Theorem 22]{Bo}.

In the formulation of the theorems below, we assume 
that $\vartheta$ ($1/2 \le \vartheta <1$) is a common level of distribution for 
primes and $E_2$-numbers. 
We then define
\begin{equation}
B = \frac2\vartheta.
\label{eq:1.18}
\end{equation}
Unconditionally, we may take $B=4$.
The Elliott--Halberstam conjecture for primes and
$E_2$-numbers is equivalent to taking $B = 2+\epsilon$.

\begin{theorem}
\label{th:1}
Let $D$ be any constant and let $L_i(x)$ $(1 \leq i \leq k)$ be
an admissible $k$-tuple of distinct linear forms.
Then there are $\nu+1$ forms among them which take simultaneously
$E_2$-numbers as values with both prime factors above $D$ 
if\ \footnote{
      For clarity, we remark that here and in subsequent results 
     (Theorems \ref{th:4}, \ref{th:5}, \ref{th:6}, 
      Corollaries \ref{cor:1}, \ref{cor:3}, \ref{cor:4}, \ref{cor:5})
      the notation $o(1)$ denotes a function $g(\nu)$ such that
      $g(\nu)\to 0$ as $\nu\to \infty$.
      }
\begin{equation}
k \geq C_1(\nu) := \frac{4 e^{-\gamma}(1 + o(1))}{B} e^{B\nu/4}.
\label{eq:1.20}
\end{equation}
\end{theorem}

\begin{theorem}
\label{th:2}
Let $\{L_1(n), L_2(n), L_3(n)\}$ be an admissible 
triplet of linear forms. 
Among these, exist two forms $L_i,L_j$ such that for infinitely
many $n$, $L_i(n), L_j(n)$ are 
both $E_2$-numbers, 
all the prime factors of which exceed $n^{1/144}$.
\end{theorem}

\begin{theorem}
\label{th:3}
Let $\{L_1(n),L_2(n)\}$ be an admissible pair of linear forms.
Then there exist infinitely many $n$ such that both $L_1(n)$
and $L_2(n)$ are $P_2$-numbers, and 
the prime factors of $L_1(n)L_2(n)$ all exceed $n^{1/10}$.
In particular, there are infinitely many integers $n$ such that
\begin{equation}
n,\, n - d \in \mathcal P_2
\label{eq:1.21}
\end{equation}
for any even integer~$d$.
\end{theorem}

Theorem~\ref{th:1} shows that, in contrast to the case of
primes, we can really prove the existence of infinitely many
blocks of $\nu$ consecutive $E_2$-numbers with a bounded
diameter (depending on $\nu$) for any given~$\nu$.

\begin{corollary} \label{cor:1}
We have for any $\nu > 0$, 
\begin{equation}
\liminf_{n \to \infty} (q_{n + \nu} - q_n) \leq
C_2(\nu) = e^{-\gamma} \nu e^{B\nu/4} (1 + o(1)).
\label{eq:1.22}
\end{equation}
\end{corollary}

Taking the admissible triplet $\{n, n+2, n+6\}$,  we see that
Theorem~\ref{th:2} implies an improvement of \eqref{eq:1.6}, namely

\begin{corollary}
\label{cor:2}
$\liminf\limits_{n \to \infty} (q_{n + 1} - q_n) \leq 6$.
\end{corollary}

The question arises: why is our present method much more
successful for $E_2$-numbers than for primes, as indicated  by
\eqref{eq:1.11}--\eqref{eq:1.12} and \eqref{eq:1.22}?

Usually sieve methods are at any rate unable to detect
$E_r$-numbers for any given $r$ due to the parity problem,  
and even in the case of $P_r$-numbers ($r$ fixed) they produce only
numbers with all prime factors larger than
\begin{equation}
X^{1/w}, \quad w > 0 \ \text{ fixed},
\label{eq:1.23}
\end{equation}
where $\text{\rm card}\, \mathcal A \sim X$, where $\mathcal A$
is a starting set containing almost primes, as emphasized by
J. H. Kan \cite{Kan1,Kan2}.
In these cases the number of almost primes considered below $N$
is $O_w(N/\log N)$ (the same as the number of primes), whereas
the true order of magnitude of the number of $P_r$-numbers (or $E_r$-numbers) is
\begin{equation}
c(r) \frac{N(\log_2 N)^{r - 1}}{(\log N)} .
\label{eq:1.24}
\end{equation}

Differently from almost all other applications of sieve methods (for
exceptions see the mentioned works \cite{Kan1}, \cite{Kan2} of
Kan), our method is able to make use of
$E_2$-numbers that satisfy
\begin{equation}
n = p_1 p_2, \quad p_1 < n^\varepsilon, 
\quad p_2 > n^{1 -\varepsilon} , 
\label{eq:1.25}
\end{equation}
for any given small $\varepsilon > 0$.
In the proof of Theorem \ref{T:S1Est}, we allow $E_2$-numbers 
with prime factors of any size.

This phenomenon 
(the {larger density of $E_2$-numbers} over primes) 
is crucial in our method.
A careful consideration of the proof of Theorem~\ref{th:1}
reveals that without taking into account the contribution of
$E_2$-numbers with \eqref{eq:1.25} for all $\varepsilon > 0$,
our method would fail to prove Theorem~\ref{th:1}.
If we exclude numbers of type \eqref{eq:1.25} for
$\varepsilon < c_0$, then we would be unable to show
Theorem~\ref{th:1}, and so Corollary~\ref{cor:1} for any $\nu >
\nu_0(c_0) \asymp c_0^{-1}$.

As we have seen in \eqref{eq:1.10}--\eqref{eq:1.11}, the level
$\vartheta$ of distribution of primes has dramatic consequences
for the strength of the result we can show about the existence
of primes in tuples.
On the other hand, the value of $\vartheta$, that is, of $B$, is
much less important in the distribution of $E_2$-numbers; only
the quantitative value $C_1(\nu)$ depends on the value of $B$,
i.e.\ of $\vartheta$.
The dependence of $C_1(\nu)$ in \eqref{eq:1.20} on $\vartheta$
is not too strong: we have in the exponent of $C_1(\nu)$
\begin{equation}
B\nu/4  = \nu/(2\vartheta) \in [\nu/2, \nu] \ \text{ for }
\vartheta \in [1/2, 1].
\label{eq:1.26}
\end{equation}

This observation has theoretical importance, for we do not
need the full strength of the Bombieri--Vinogradov theorem.
Moreover,
it can be used to generalize the results of Theorem~\ref{th:1}
for a situation when 
$\vartheta = \varepsilon \Leftrightarrow B= 2/\varepsilon$, for example.
We remark that, contrary to this, the proof of
\eqref{eq:1.1} would break down if we had  just a fixed
$\vartheta < 1/2$ at our disposal, even if this value were
very close to~$1/2$.

The case of $\vartheta$ being small occurs when we would like to
find blocks of bounded length of $E_2$-numbers in short
intervals of type
\begin{equation}
\big[N, N + N^{7/12 + \varepsilon}\big], \quad \varepsilon > 0 \
\text{ fixed.}
\label{eq:1.27}
\end{equation}

In this case it was proved by Perelli, Pintz and Salerno
\cite{PPS} in 1985 that one has a short interval version of
Bombieri--Vinogradov theorem for intervals of type
\eqref{eq:1.27} where (surprisingly) we can choose $\vartheta$
as a fixed positive constant, $\vartheta = 1/40$ for any
$\varepsilon > 0$.
This was improved two years later by Timofeev \cite{Tim} to
\begin{equation}
\vartheta = 1/30 \Longleftrightarrow B = 60.
\label{eq:1.28}
\end{equation}

The result proved in \cite{PPS} reads as
\begin{equation}
\sum_{q \leq Q} \max_{(a,q) = 1} \max_{h \leq y} \max_{x/2 < z
\leq x} \bigg| \sum_{\substack{p \equiv a(\text{\rm mod}\, q)\\
z < p \leq z + h}} \log p - \frac{h}{\varphi(q)} \bigg| \ll
\frac{y}{(\log x)^A} , 
\label{eq:1.29}
\end{equation}
where $y = x^{7/12 + \varepsilon}, \quad
Q = x^\vartheta(\log x)^{-D}, \quad A \text{ arbitrary,} \quad D
= D(A)$.

The method of proof of both works \cite{PPS} and \cite{Tim} uses
Heath-Brown's identity, therefore the analogue of \eqref{eq:1.29}
can be proved {\it mutatis mutandis} for $E_2$-numbers as well.
Accordingly, we will prove

\begin{theorem}
\label{th:4}
Under the conditions of Theorem~\ref{th:1} we can find $\nu+1$
linear forms $L_{i_j}(n)$ which take $E_2$-numbers with both
prime factors above $D$ and for any $\varepsilon > 0$ we can
require
\begin{equation}
n \in \bigl[N, N + N^{7/12 + \varepsilon}\bigr]
\label{eq:1.30}
\end{equation}
if
\begin{equation}
k \geq C_3(\nu) := \frac{ e^{-\gamma}(1 + o(1))}{15} e^{15\nu},
\quad N > N_0(k, \varepsilon, D).
\label{eq:1.31}
\end{equation}
\end{theorem}

\begin{corollary}
\label{cor:3}
Let $\nu \in \mathbb Z^+$, $\varepsilon > 0$.
If $N > N_1(\nu, \varepsilon)$, then there exist $n, m \in
\mathbb Z^+$ such that
\begin{equation}
n \in \bigl[N, N + N^{7/12 + \varepsilon}\bigr]
\label{eq:1.32}
\end{equation}
\begin{equation}
n < q_m < q_{m + 1} < \dots < q_{m + \nu} < n +  e^{-\gamma}\nu e^{15\nu}(1 + o(1)) .
\label{eq:1.33}
\end{equation}
\end{corollary}

We can further restrict our $E_2$-numbers $p_1 p_2$ to be of the
form
\begin{equation}
p_1 p_2 = x^2 + y^2\quad (x,y \in \mathbb Z) \Longleftrightarrow
p_1, p_2 \equiv 1\ (\text{\rm mod}\, 4).
\label{eq:1.34}
\end{equation}

It is relatively easy to show the existence of infinitely many
families of triplets of consecutive integers 
that are sums of two squares--see \cite{CD} or \cite{Ho} for a more general 
result.
A modification of our proof of Theorem~\ref{th:1} shows
that Theorem~\ref{th:1} remains essentially valid 
for $E_2$-numbers which are sums of two squares.

\begin{theorem}
\label{th:5}
Under the conditions of Theorem~\ref{th:1}, we have infinitely
many $n$ such that at least $\nu+1$ linear forms $L_{i_j}(n)$
take simultaneously $E_2$-values (with both prime factors
above~$D$) which are sums of two squares, if
\begin{equation}
k \geq C_4(\nu) := \frac{4e^{-\gamma}(1 + o(1))}{B } e^{B\nu}.
\label{eq:1.35}
\end{equation}
\end{theorem}

\begin{corollary}
\label{cor:4}
Under the conditions of Theorem~\ref{th:5} we have at least
$\nu$ $E_2$-numbers which are sums of two squares infinitely
often in intervals of length $K$ if
\begin{equation}
K \geq C_5(\nu) :=  {4e^{-\gamma}\nu} e^{B\nu} (1 + o(1)).
\label{eq:1.36}
\end{equation}
\end{corollary}

Finally we can combine the results of Theorems~\ref{th:4} and
\ref{th:5} to have $\nu$ $E_2$-numbers which are sums of two
squares if the value $n$ is localized in a short interval of
type~\eqref{eq:1.30}.

\begin{theorem}
\label{th:6}
In Theorem~\ref{th:4} we may require that the $\nu$ $E_2$-values
of the linear forms should be sums of two squares if we have in
place of \eqref{eq:1.31} the restriction
\begin{equation}
k \geq C_6(\nu) := \frac{e^{-\gamma}(1 + o(1))}{15} e^{60\nu}. 
\label{eq:1.37}
\end{equation}
\end{theorem}

\begin{corollary}
\label{cor:5}
Let $\{q'_n\}^\infty_1$ denote the set of $E_2$ numbers which
can be written as sums of two squares.
Then Corollary~\ref{cor:3} is true if \eqref{eq:1.33} is
replaced by
\begin{equation}
n < q'_m < q'_{m+1} < q'_{m+\nu} < n + 
{4e^{-\gamma}\nu} e^{60\nu}(1 + o(1)).
\label{eq:1.38}
\end{equation}
\end{corollary}

It will be clear from the proofs that all of the above theorems and corollaries
remain true if we require that all of the constructed $E_2$-numbers have
both of their prime factors exceeding some specific constant.
Indeed, this holds more generally if both prime factors exceed some $Y(N)$
with $\log Y(N)/\log N\to 0$ as $N\to \infty$. 

Our methods open the way towards a new, simpler and unified
treatment of some conjectures of Erd\H os \cite{Er}  on consecutive
integers, the most well-known of them being the Erd\H os--Mirsky
\cite{EM}
conjecture, which states that
\begin{equation}
d(n) = d(n + 1) \ \text{ infinitely often (i.o.)};
\label{eq:1.39}
\end{equation}
the two others being the analogous conjectures with $d(n)$
replaced by the functions $f(n) = \omega(n)$ or $\Omega(n)$
(number of prime divisors of $n$ counted without and with
multiplicity, respectively). Similar to \eqref{eq:1.3}--\eqref{eq:1.4},
the parity problem seemed to prevent a solution of these
conjectures. However, as observed by Spiro \cite{Spi} and
Heath-Brown \cite{HB2}, the parity problem can be ``sidestepped,''
and it is possible to prove the conjectures without specifying the
common value of $f(n) = f(n + 1)$ (or even its {parity}) for
the relevant functions $f = d$ or $\Omega$. 
Recently, the same was shown for 
$f = \omega$   by Schlage-Puchta \cite{SP}.

In the next part of this series we will show these conjectures
in the stronger form, where we can specify the common value
$f(n) = f(n + 1)$, even in a nearly arbitrary way in case of
$\omega$ or $\Omega$, thereby overcoming the parity problem in
these cases.
We will prove

(i) $d(n) = d(n + 1) = A$ holds i.o.\ for any $A$ with $24 \mid A$,

(ii) $\omega(n) = \omega(n + 1) = A$ holds i.o.\ for any $A \geq
3$,

(iii) $\Omega(n) = \Omega(n + 1) = A$ holds i.o.\ for any $A
\geq 4$.

Further we can show the analogous statements in case of an
arbitrary shift $b$ in place of $1$, if $f = \omega$ or $\Omega$
(where the lower bound for $A$ may depend on $b$ in case of
$\omega$, and should be replaced uniformly by $5$ in case of
$\Omega$) and for every shift $b \not\equiv 15$ $(\text{\rm
mod}\, 30)$ for the divisor function.

This generalization was proved for every $b$ by Pinner
\cite{Pin} in 1997 (without specifying the common value of $f$)
for $f = d$ and $\Omega$; however, the method used by Schlage-Puchta
for $\omega$ does not work for general $b$.
On the other hand, Buttkewitz \cite{Bu} has recently proved that an
analogous result holds for an infinite set $\mathcal{B}$ of possible shifts $b$.

\section{Notation and Preliminary Lemmas} \label{S:Lemmas}

Most of our notation will be introduced as needed, but it is useful to make the following comments here.
Throughout this paper, we use $k$ to denote an integer $k\ge 2$,
$\cL$ to denote an admissible $k$-tuple of linear forms, and
$P$ to denote a polynomial. The constants implied by $``O''$ and $``\ll''$
may depend on $k,\cL,$ and $P$. 
$\tau_k(n)$ denotes the number of ways of writing $n$ as product of $k$ factors.
$\omega(n)$ is the number of distinct prime factors of $n$. 
$\phi(n)$ and $\mu(n)$ are the usual functions of Euler and M\"obius, respectively.
The letters $N$ and $R$ denote real numbers regarded as tending to infinity, 
and we always assume that $R\le N^{1/2}$. 

To count $E_2$-numbers, we introduce the following function $\beta$. 
Let $Y$ be a real number with  $1\le Y \le N^{1/4}$, 
and  define
\begin{equation} \label{eq:betaDef}
\beta(n)= 
\begin{cases} 
 1 & \text{ if $n=p_1p_2$, $Y< p_1 \le N^{1/2} < p_2$,}\\
 0 & \text{ otherwise.}
 \end{cases}
 \end{equation}

The notation $\pi(x)$ is commonly used to denote the number of primes up to $x$,
and $\pi(x;q,a)$ denotes the number of primes up to $x$ that are congruent to $a \pmod q$. 
For our purposes, it is convenient to define the following related quantities. 
\begin{align*} 
 \pi^\flat (x) &  =  \sum_{x< p \le 2x}  1 = \pi(2x)-\pi(x) \\
 \pi^\flat (x;q,a) &  = \sum_{\substack{ x < p \le 2x \\ p\equiv a \pmod q}} 1= \pi(2x;q,a)-\pi(x;q,a) \\
\pi_\beta(x) & = \sum_{x < n \le 2x} \beta(n)  \\ 
\pi_{\beta,u}(x) & = \sum_{\substack{ x < n \le 2x \\ (n,u)=1}} \beta(n) \\
 \pi_\beta(x;q,a) & = \sum_{\substack{ x < n \le 2x \\ n \equiv a \pmod q}} \beta(n) \\
\end{align*}

As mentioned in the introduction, we will employ results on the level of distribution for both
prime numbers and $E_2$-numbers. For primes, define
\begin{equation*}
\Delta(x;q,a) = \pi^\flat(x;q,a) - \frac{1}{\phi(q)} \pi^\flat(x)
\end{equation*}
and 
\begin{equation*}
\Delta^*(x;q) = \max_{y\le x} \max_{a; (a,q)=1} |\Delta(y;q,a)|
\end{equation*}

\begin{lemma} \label{L:Lemma3}
Assume that the primes have level of distribution $\vartheta$, $\vartheta \le 1$.
For every $A>0$ and for every fixed integer $h\ge 0$,
 there exists $C=C(A,h)$ such that 
 if $Q\le x^{\vartheta} (\log x)^{-C}$, then
\begin{equation*}
\sum_{q\le Q}  \mu^2(q) h^{\omega(q)}  \Delta^*(x;q)  \ll_A x (\log x)^{-A}.
\end{equation*}
\end{lemma}

By the Bombieri-Vinogradov Theorem, this lemma is unconditional for $\vartheta \le 1/2$.
The incorporation of the factor $h^{\omega(q)}$ is familiar feature in sieve applications;
see \cite[Lemma 3.5]{HR}, for example.

For the function $\beta$, we define 
 \begin{align*}
 \Delta_\beta(x;q,a) = &
 \sum_{\substack{  x < n \le 2x \\ n \equiv a \pmod q}} \beta(n) -
   \frac{1}{\phi(q)} \sum_{\substack{ x < n \le 2x \\ (n,q)=1} } \beta(n) =
 \pi_\beta(x;q,a)- \frac{1}{\phi(q)} \pi_{\beta,q} (x), \\
 \Delta^*_\beta(x;q) = & \max_{y\le x} \max_{a; (a,q)=1} |\Delta_\beta (y;q,a)|.
\end{align*}

\begin{lemma} \label{L:Lemma4}
Assume that $E_2$-numbers have a level of distribution $\vartheta$, 
$\vartheta \le 1$.
For every $A>0$ and  for every fixed integer $h\ge 0$,
 there exists $C=C(A,h)$ such that  if  
 $Q\le x^{\vartheta} (\log x)^{-C}$, then
\begin{equation} \label{E:BVE2}
\sum_{q\le Q}  \mu^2(q) h^{\omega(q)} \Delta^*_\beta(x;q)   \ll_A x (\log x)^{-A}.
\end{equation}
\end{lemma}

This follows from a general result of Motohashi \cite{Mot} when $\vartheta \le 1/2$. 
When $1/2 < \vartheta \le 1$, Lemmas \ref{L:Lemma3} and \ref{L:Lemma4}
are both hypothetical.

Our next lemma is central to the estimation of the sums that arise
in Selberg's sieve.

\begin{lemma} \label{L:Wirsing}
Suppose that $\gamma$ is a multiplicative function,
and suppose that there are positive real numbers 
$\kappa,A_1,A_2, L$ such that
\begin{equation} \label{E:Omega1}
	0\le \frac{\gamma(p)}{p} \le 1-\frac{1}{A_1},
\end{equation}
and 
\begin{equation} \label{E:Omega2}
	-L \le \sum_{w\le p < z} 
	         \frac{\gamma(p) \log p}{p}
		 -\kappa \log \frac{z}{w} 
		\le A_2
\end{equation}
if $2 \le w \le z$. 
Let $g$ be the multiplicative function defined by 
\begin{equation} \label{E:gDefinition}
    g(d) = \prod_{p|d} \frac{\gamma(p)}{p-\gamma(p)}.
\end{equation}
Then 
\begin{equation*}
    \sum_{d < z} \mu^2(d) g(d) =
     c_{\gamma} \frac{(\log z)^{\kappa}}{\Gamma(\kappa+1)}
        \left\{ 1+O\left( \frac{L}{\log z} \right) \right\},
\end{equation*} 
where
\begin{equation*}
    c_{\gamma} = 
     \prod_p \left( 1-\frac{\gamma(p)}{p} \right)^{-1}
             \left( 1-\frac{1}{p} \right)^{\kappa}.
\end{equation*}
\end{lemma}

This is a combination of Lemmas 5.3 and 5.4 of Halberstam and 
Richert's book \cite{HR}. 
In \cite{HR}, the hypothesis \eqref{E:Omega1} is denoted
$(\Omega_1)$, and hypothesis \eqref{E:Omega2} is denoted
$(\Omega_2(\kappa,L))$.
The constant implied by ``$O$'' may depend on 
$A_1,A_2,\kappa$, but it is independent of $L$. 

Our next lemma is a variant of the previous one with the terms $g(d)$ weighted by an 
appropriate function.

\begin{lemma} \label{L:WirsingT}
Assume the hypotheses of Lemma \ref{L:Wirsing}.
Assume also that $F:[0,1]\to \bR$ is a piecewise differentiable function.
Then
\begin{align} \label{E:WirsingT}
\sum_{d< z } \mu^2(d) g(d) F\left( \frac{\log z/d}{\log z} \right)= 
 c_\gamma &  \frac{(\log z)^\kappa}{\Gamma(\kappa)}   \int_0^1 F(1-x) x^{\kappa-1} dx \\
  &    + O \left(c_\gamma L M(F) (\log z)^{\kappa-1}\right), \notag
 \end{align}
 where  $ M(F) = \sup \{ (|F(x)|+|F'(x)|): 0\le x \le 1\}.$
 The constant implied by ``O'' may depend on $A_1,A_2, $ and $\kappa,$  but it is 
 independent of $L$ and $F$.
 \end{lemma}
 
  \begin{proof}
The left-hand side of the proposed conclusion is
 \begin{equation*}
  \int_{1^-}^z F\left(\frac{\log z/u}{\log z}\right) dG(u),
\end{equation*}
where  
\begin{equation*}
  G(u) = \sum_{d< u} \mu^2(d) g(d) = c_\gamma \frac{(\log u)^\kappa}{\Gamma(\kappa+1)}+ E(u), \end{equation*}
and $E(u) \ll c_\gamma L (\log 2u)^{\kappa-1}$ by the previous lemma.
Consequently, the sum in \eqref{E:WirsingT} may be written as 
\begin{equation*}
 \int_1^z  F\left(\frac{\log z/u}{\log z}\right) d c_\gamma \frac{(\log u)^\kappa}{\Gamma(\kappa+1)} +
 \int_{1^-}^z  F\left(\frac{\log z/u}{\log z}\right) dE(u).
\end{equation*}
In the first integral, we make the change of variables $u=z^x$; this gives the first term.
We use integration by parts on the second integral to obtain
\begin{equation*}
  \int_{1^-}^z  F\left(\frac{\log z/u}{\log z}\right) dE(u)=
  \left.F\left(\frac{\log z/u}{\log z}\right) E(u) \right]_{1^-}^z +
  \int_1^z E(u) F'\left(\frac{\log z/u}{\log z}\right) \frac{du}{u \log z}.
\end{equation*}
The desired result now follows by using the above-mentioned bound for $E(u)$.
 \end{proof}

\section{Initial Considerations}\label{sec:2}  

Let $\cL=\{L_1,L_2, \ldots, L_k\}$ be an admissible $k$-tuple of 
linear forms; i.e., a set of linear forms satisfying the 
conditions of \eqref{eq:1.15} and \eqref{eq:1.16}. 
Admissibility can also be defined in terms of solutions of 
congruences. Define
\begin{equation*} 
    P_{\cL}(n)=\prod_{i=1}^{k} L_{i}(n)  = (a_1n+b_1) \ldots (a_k n+b_k),
\end{equation*}
and for each prime $p$ define
\begin{align*}
    \Omega_p(\cL) & =
    \{ n: 1\le n \le p \text{ and } P_{\cL}(n)\equiv 0 \pmod p\},\\
    \nu_{p}(\cL) & = \#\Omega_{p}(\cL).
\end{align*}
The condition that $\cL$ is admissible is equivalent to requiring that
\begin{equation} \label{eq:2.0}
    \nu_{p}(\cL) < p
\end{equation}
for all primes $p$.
We always have $\nu_{p}(\cL) \le k$, so the above condition is 
automatic for any prime $p>k$.

The singular series connected to $\cL$ is defined as 
\begin{equation*}
    \fS(\cL) = 
    \prod_{p}
      \left( 1 - \frac{\nu_{p}(\cL)}{p} \right)
      \left( 1 - \frac{1}{p} \right)^{-k}.
\end{equation*}
The product converges because
$\nu_p =k$ for all but
finitely many primes $p$, and 
$\cL$ is admissible if and only if $\fS(\cL)\ne 0$. 

Next, we examine those  primes $p$ satisfying $\nu_{p} <k$. 
First of all, if $p|a_i$ for some $i$, then 
the congruence $a_{i}n+b_i \equiv 0 \pmod p$ will have no solutions,
and so $\nu_p < k$ in this case. 
Now suppose that $p\nmid a_i$ for all $i$. In this instance,
 $\nu_{p}< k$ if
and only if there are indices $i,j \, (i\ne j)$ such that
\begin{equation*}
    \overline{a}_i b_i \equiv \overline{a}_j b_j \pmod p,
\end{equation*} 
where $\overline{a}_i$ denotes the multiplicative inverse of $a_i \pmod p$.
We therefore see that $\nu_{p}<k$ if and only if $p|A$,
where
\begin{equation} \label{eq:2.1}
    A=A(\cL)= \prod_{i=1}^{k}  a_{i}  
                  \prod_{1\le i < j \le k} |a_ib_j-a_jb_i|.
\end{equation}

For technical reasons, it is useful to adopt the normalization 
introduced by Heath-Brown\cite{HB}.
For each prime $p|A$, there is an integer $n_p$ such that 
$p\nmid P_{\cL}(n_p)$. By the Chinese Remainder Theorem,
there is an integer $B$ such that $B\equiv n_p \pmod p$ for all 
$p|A$. For  $i=1,\ldots,k$, we define
\begin{equation*}
L_i' (n) = L_i(An+B)=  a'_in + b'_i,
\end{equation*}
where $a'_i = a_i A $ and $b'_i= L_i(B) = a_i B+ b_i$.
Set
\begin{equation*}
\cL'=\{ L_1', \ldots , L_k' \}.    
\end{equation*}
We claim that
\begin{equation} \label{eq:2.2}
    \nu_p(\cL')=
    \begin{cases}
	k & \text{ if $p\nmid A$, and}\\
	0 & \text{ if $p \mid A$.}
    \end{cases}
\end{equation}
To justify this claim, we assume first that $p|A$. 
Then 
$$
L'_i(n)\equiv L_i(B) \equiv L_i(n_p) \not\equiv 0 \pmod p
$$
for all integers $n$, and so $\nu_p(\cL')=0$. 
Next assume that $p\nmid A$. As noted before,
$\nu_p(\cL')<k$ if and only if $p|(a'_ib'_j - a'_jb'_i)$ for some 
choice of $i,j$ with $1\le i < j \le k$. 
However,
\begin{equation*}
    a'_ib'_j - a'_jb'_i=
    \det
     \left| 
         \begin{matrix} a_i A & a_j A \\ a_iB+b_i & a_jB+b_j \end{matrix}
      \right| =
      A \det
       \left| 
         \begin{matrix} a_i  & a_j  \\ b_i & b_j \end{matrix}
      \right|,
\end{equation*}
and this is not divisible by $p$.

For brevity, it is useful to relabel $L'_{i}$ as $L_i$ and to assume 
the following hypothesis.

\begin{HypoA} 
$\cL=\{L_1, \ldots, L_k\}$ is an admissible $k$-tuple of linear forms. 
    The functions $L_i(n)=a_i n+b_i (1\le i\le k)$ have 
integer coefficients with $a_i>0$. Each of the coefficients $a_i$ is 
composed of the same primes, none of which divides any of the $b_i$. 
If $i\ne j$, then any prime factor of $a_ib_j-a_jb_i$ divides each of 
the $a_i$.
\end{HypoA}

For sets of linear forms $\cL$ satisfying Hypothesis A, we re-define
$$
A=\prod_i a_i.
$$
In this case, 
\begin{equation*}
    \fS(\cL) = \prod_{p|A} \left( 1-\frac{1}{p} \right)^{-k}
               \prod_{p\nmid A} \left( 1-\frac{k}{p} \right)
	                   \left( 1-\frac{1}{p} \right)^{-k}.
\end{equation*}
Note that by \eqref{eq:2.0} and \eqref{eq:2.2},  
\begin{equation} \label{eq:2.2a}
p\le k \implies p|A,
\end{equation}
so  $\fS(\cL)$ is indeed positive. 
  
The primary tool for proving Theorems 1 through 3 is an adaptation of 
the basic construction of Goldston, Pintz, and \Yildirim.
Let $\cL=\{L_1,L_2, \ldots , L_k\}$ be a set of linear forms 
satisfying Hypothesis A, and let $\beta$ be as defined in \eqref{eq:betaDef}. 
For the proof of Theorem 1, we seek an asymptotic
formula for the sum 
\begin{equation} \label{eq:cSDef}
\cS = \sum_{N< n \le 2N} 
 \left\{ \sum_{j=1}^k \beta(L_j(n)) \,\,-\,\,\nu \right\} 
    \left(\sum_{d|P_{\cL}(n)} \lambda_d \right)^2,
\end{equation}
where the $\lambda_d$'s are real numbers to be chosen in due course.
The significance of $\cS$ is that  a value of $n$ contributes a 
positive amount only if at least $\nu+1$ elements 
of the set $\{L_1(n), \ldots, L_k(n)\}$ are $E_2$-numbers.

We  immediately decompose $\cS$ as
\begin{equation} \label{eq:SDecomp}
\cS = \sum_{j=1}^k \cS_{1,j} - \nu \cS_0,
\end{equation}
where
\begin{equation*}
\cS_{1,j} = \sum_{N< n \le 2N} \beta(L_j(n)) 
 \left(\sum_{d|P_{\cL}(n)} \lambda_d \right)^2,
\end{equation*}
and 
\begin{equation*}
\cS_{0} = \sum_{N< n \le 2N} 
 \left(\sum_{d|P_{\cL}(n)} \lambda_d \right)^2,
\end{equation*}

The motivation for the use of the coefficient $\lambda_d$
comes from the realm of the Selberg sieve. More specifically, consider 
the problem of bounding the number of
$n$ for which all of $L_1(n), \ldots, L_k(n)$ are prime. 
Start from the observation that if $\lambda_1=1$ 
and $\lambda_d=0$ for $d>N$, then 
\begin{align}
\sum_{\substack{N < n \le 2N \\ \text{ all $L_i(n)$ prime} }} 1 
\notag
 & \le
\sum_{N < n \le 2N }\left(\sum_{d|P_{\cL}(n)} \lambda_d \right)^2 \\
& = \sum_{d,e} \lambda_d \lambda_e 
      \sum_{\substack{ N < n \le 2N \\ [d,e]|P_{\cL}(n) }} 1. 
      \label{eq:2.3}
\end{align}

As we noted in the introduction, we take 
$\nu_p(\cL)$ to be 
the number of solutions of $P_{\cL}(n) \equiv 0 \pmod p$. 
We extend this definition to arbitrary squarefree $d$ by multiplicativity. 
Consequently,
\begin{equation*}
\sum_{\substack{ N < n \le 2N \\ d|P_{\cL}(n)}} 1=
N \frac{\nu_d(\cL)}{d} + O\left(k^{\omega(d)}\right)
\end{equation*}
for squarefree $d$.
Returning to \eqref{eq:2.3}, we find that the expression there is
\begin{equation*}
N \sum_{d,e} \frac{\lambda_d \lambda_e \nu_{[d,e]}(\cL)}{[d,e]} +
O\left(\sum_{d,e} |\lambda_d \lambda_e| k^{\omega([d,e])} \right).
\end{equation*}
We control the size of the error term by specifying that 
$\lambda_d=0$ if $d\ge R$, where $R$ will be chosen later. 
Moreover, the terms with $([d,e],A)>1$ make no contribution since
$\nu_{[d,e]}(\cL)=0$ for these terms. Accordingly, we restrict the sum to 
terms with $(d,A)=(e,A)=1$. 
It is also convenient to specify that
\begin{equation*}
\lambda_d = 0 \text{ if $d$ is not squarefree.}
\end{equation*}

The coefficient of $N$ in the main term may be rewritten as
\begin{equation} \label{E:BilinearForm1}
    \sumprime_{d,e} \frac{\lambda_d \lambda_e}{f([d,e])},
\end{equation}
where $\sumprime_\null$ denotes that the sum is over all values of 
the indices that are relatively prime to $A$, and 
\begin{equation} \label{E:fDef}
f(d) = \frac{d}{\nu_d(\cL)} = \frac{d}{\tau_k(d)}=  \prod_{p|d} \frac{p}{k}
\end{equation}
for squarefree $d$ with $(d,A)=1$. 

The typical approach in the Selberg sieve 
is to choose the $\lambda_d$ to minimize the form
in \eqref{E:BilinearForm1}. To make this problem feasible, one needs 
to diagonalize this bilinear form. This can be done 
by making a change of variables
\begin{equation} \label{E:yrDefinition1}
    y_r= \mu(r) f_1(r) \sumprime_{d} \frac{\lambda_{dr}}{f(dr)},
\end{equation}
where $f_1$ is the multiplicative function defined by $f_1=f*\mu$.
In other words, 
\begin{equation}\label{E:f1Def}
f_1(d) = \prod_{p|d} \frac{p-k}{k}
\end{equation}
whenever $d$ is squarefree and $(d,A)=1$. 
(Note that the sum in \eqref{E:yrDefinition1} is finite because
$\lambda_d=0$ for $d>R$. 
Note also that there is an implicit
condition $(d,r)=1$ because $\lambda_{dr}=0$ if
$dr$ is not squarefree.) The sum in \eqref{E:BilinearForm1} is 
then transformed into
\begin{equation*}
    \sumprime_{r} \frac{y_r^2}{f_1(r)},
\end{equation*}
and the bilinear form is minimized by taking
\begin{equation} \label{E:Sy}
    y_r = \mu^2(r) \frac{\lambda_1}{V}
\end{equation}
when $r<R$ and $(r,A)=1$, where
\begin{equation*}
    V=\sumprime_{r < R} \frac{\mu^2(r)}{f_1(r)}.
\end{equation*}
The minimum of the form in \eqref{E:BilinearForm1} is then seen
to be
\begin{equation*}
\frac{\lambda_1^2}{V}.
\end{equation*}
One usually assumes that $\lambda_1=1$, but this is not an essential
element of the Selberg sieve, and it is sometimes useful to assign 
some other nonzero value to $\lambda_1$.

Using M\"obius inversion, one can easily show that
\begin{equation} \label{eq:YToLambda}
\lambda_d = \mu(d) f(d) \sumprime_{r} \frac{y_{rd}}{f_1(rd)}.
\end{equation}
Consequently, specifying a choice for $\lambda_d$ is equivalent to 
specifying a choice for $y_r$. 
Our choice of $\lambda_d$ is different from the choice implied by (3.12),
and it is most easily described in terms of $y_r$.
We will take
\begin{equation} \label{E:yrChoice}
    y_r=
     \begin{cases}
	 \displaystyle  \mu^2(r) \fS(\cL) 
                        P\left( \frac{\log R/r}{\log R} \right) &
	  \text{ if $r< R$ and $(r,A)=1$,}\\
	  0 &\text{ otherwise.}
      \end{cases}
\end{equation}
Here, $P$ is a polynomial to be determined later.

Our estimate for $\cS$ follows from the following two results.

\begin{theorem} \label{T:S0Est}
Suppose that $\cL$ is a set of linear forms satisfying Hypothesis A.
Suppose that  $\lambda_d$ is given by   \eqref{eq:YToLambda}
and \eqref{E:yrChoice}. 
There is some constant $C$ such that if $R\le N^{1/2} (\log N)^{-C}$, then
\begin{equation*}
\cS_0 =  \frac{\fS(\cL) N (\log R)^k}{(k-1)!} J_0 +
    O\left( N(\log N)^{k-1} \right)
\end{equation*}
where
\begin{equation*}
J_0=\int_0^1 P(1-x)^2 x^{k-1} dx.
\end{equation*}
\end{theorem}

\begin{theorem} \label{T:S1Est}
Suppose that the primes and the $E_2$-numbers have a common level of distribution
$\vartheta \le 1$, and let $\cL$ be a set of linear forms satisfying Hypothesis A.
Suppose that $\lambda_d$ is given by \eqref{eq:YToLambda} and \eqref{E:yrChoice}, 
and let the polynomial $\Pt$ is defined as 
\begin{equation*}
\Pt(x)=\int_0^x P(t) dt.
\end{equation*}
There is some constant $C$ such that if $R=N^{\vartheta/2} (\log N)^{-C}$, then
\begin{equation*}
\cS_{1,j} = \frac{\fS(\cL) N (\log R)^{k+1}}{(k-2)!( \log N)} (J_1+J_2+J_3)
 + O\left( N (\log \log N) (\log N)^{k-1} \right),
\end{equation*}
where $Y= N^\eta$, $B=2/\vartheta$, and 
\begin{align*}
J_1 = & \int_{B\eta}^1 \frac{B}{y(B-y)} 
       \int_0^{1-y} \left( \Pt(1-x)-\Pt(1-x-y) \right)^2 x^{k-2} 
          dx \, dy,\\
J_2 = & \int_{B\eta}^1 \frac{B}{y(B-y)}
       \int_{1-y}^1 \Pt(1-x)^2 x^{k-2} dx \, dy,\\
J_3 = & \int_1^{B/2} \frac{B}{y(B-y)} \int_0^1 \Pt(1-x)^2 x^{k-2} dx \, dy.
\end{align*}
\end{theorem}

Finally, we mention the following result, which is needed for the proof of 
Theorem \ref{th:3}.

\begin{theorem} \label{th:9}
Assume the hypotheses of Theorem \ref{T:S1Est}. 
Let $\varpi$ denote the characterstic function of the primes; i.e.,
$\varpi(p)=1$ if $p$ is a prime and $\varpi(p)=0$ otherwise. 
There is some constant $C$ such that if $R\le N^{\vartheta/2} (\log N)^{-C}$, then
\begin{align*}
\sum_{N< n \le 2N} \varpi(L_j(n))  \left(\sum_{d|P_{\cL}(n)} \lambda_d \right)^2 
 =&  \frac{\fS(\cL) N (\log R)^{k+1}}{(k-2)! (\log N)} J_{\varpi} \\
 & \text{\quad\quad} + O\left( N   ( \log\log N) (\log N)^{k-1}  \right),\\
\end{align*}
where
\begin{equation*}
J_{\varpi}= \int_0^1 \Pt(1-x)^2 x^{k-2} dx.
\end{equation*}
\end{theorem}

This result is very similar to Theorem 1.6 of \cite{GGPY1} and to Theorem \ref{T:S1Est};
consequently, we will
give just a short sketch of the proof in Section \ref{sec:7}.

\section{Proof of Theorem \ref{T:S0Est}} \label{sec:3}

From the definition of $\cS_0$, we see  that
\begin{align} 
\cS_0= \sumprime_{d,e} \lambda_{d} \lambda_{e}
       \sum_{\substack{N < n \le 2N \\ [d,e]|P_{\cL}(n) }} 1  
& = N \sumprime_{d,e} \frac{\lambda_{d} \lambda_{e}}{f([d,e])} 
  + O\left(
       \sumprime_{d,e} |\lambda_{d}\lambda_{e} r_{[d,e]}|\right)  \notag \\
& = NS_{01} + O(S_{02}),\label{E:S0Decomp}
\end{align} 
say, where
\begin{equation} \label{E:rdDef}
    r_d= \sum_{\substack{ N< n\le 2N \\ d|P_{\cL}(n)}} 1 - \frac{N}{f(d)}.
\end{equation}
Now 
\begin{equation*}
S_{01} =  \sumprime_{d,e} \frac{\lambda_d\lambda_e}{f(d)f(e)} 
        \sumprime_{\substack{r|d\\r|e}} f_1(r) 
    =  \sumprime_r f_1(r) 
          \left(\sumprime_d\frac{\lambda_{dr}}{f(dr)}\right)^2 
    =  \sumprime_r \frac{\mu^2(r) y_r^2}{f_1(r)}.
\end{equation*}

We use Lemma \ref{L:WirsingT} with 
\begin{equation*}
\gamma(p)= 
  \begin{cases}
     k & \text{ if $p\nmid A$,} \\
     0 & \text{ if $p| A$}
  \end{cases}
\end{equation*}
and $\kappa=k, F(x)=P(x)^2$. 
We deduce that 
\begin{equation*}
S_{01} = \frac{\fS(\cL) (\log R)^k}{(k-1)!} 
                      \int_0^1 P(1-x)^2 x^{k-1} dx 
                +O\left( (\log R)^{k-1} \right).
\end{equation*}

For $S_{02}$,  we first note the bound
\begin{equation*}
|r_{[d,e]} | \le k^{\omega([d,e])}
\end{equation*}
that follows from \eqref{E:rdDef}. We will later establish the bound
\begin{equation} \label{eq:lambdabound}
|\lambda_d|\ll (\log R)^k
\end{equation}
whenever $d\le R$ and $d$ is squarefree.
Assuming this momentarily, we find that 
\begin{align*}
S_{02}&  \ll (\log R)^{2k} \sum_{d,e< R} \mu^2(d) \mu^2(e) k^{\omega([d,e])} 
            \ll (\log R)^{2k} \sum_{r < R^2} \mu^2(r) (3k)^{\omega(r)}  \\
         & \ll R^2 (\log R)^{2k}\sum_{r< R^2} \frac{\mu^2(r) (3k)^{\omega(r)}}{r} 
            \ll R^2 (\log R)^{2k}\prod_{p<R^2} \left(1 + \frac{3k}{p} \right) \\
          &\ll R^2 (\log R)^{5k}.
\end{align*}
Therefore $S_{02}\ll N$ if $R\le N^{1/2} (\log N)^{-3k}$. 

To finish, we need to establish the bound \eqref{eq:lambdabound} on $|\lambda_d|$.
From \eqref{eq:YToLambda}, we see that if $d\le R$ and $d$ is squarefree, then 
\begin{align*}  
 |\lambda_d|  = &\fS(\cL)  \frac{f(d)}{f_1(d)}  \sumprime_{\substack{r < R/d\\ (r,d)=1}} 
                                              \frac{\mu^2(r)}{f_1(r)} P\left(\frac{\log R/dr}{\log R}\right) 
                                              \\
                      \le  & \fS(\cL) \sup_{0\le u \le 1} |P(u)|
                                      \sum_{\delta| d} \frac{1}{f_1(\delta)}
                                          \sumprime_{\substack{r < R/\delta\\ (r,d)=1}} 
                                                \frac{\mu^2(r)}{f_1(r)}  \\
                       \ll  & \sumprime_{r<R} \frac{\mu^2(r)}{f_1(r)} \ll (\log R)^k,  
\end{align*}
where we have used Lemma \ref{L:Wirsing}  in the last line.

\section{Proof of Theorem \ref{T:S1Est}--Initial Steps} \label{sec:4}

From the definition of $\cS_{1,j}$, we see  that
\begin{equation} \label{eq:50}
\cS_{1,j} = \sum_{d,e} \lambda_d \lambda_e 
   \sum_{\substack{ N < n \le 2N \\ [d,e]|P_{\cL}(n) }} \beta(L_j(n)).
\end{equation}
We analyze the inner sum in the next lemma.

\begin{lemma} \label{L:Lemma5}
Suppose that $u$ is squarefree, $(u,A)=1$, and all prime divisors of $u$ are less than $R$.
Define 
\begin{equation*}
M_j(u) = \sum_{\substack{N< n \le 2N \\ u|P_\cL(n)}} \beta(L_j(n)).
\end{equation*}
Then
 \begin{align}
M_j(u) = & 
 \frac{ \tau_{k-1}(u) }{\phi(a_ju)} 
     \sum_{\substack{ Y < p \le N^{1/2} \\ p\nmid u}} \pi^\flat(a_j N/p) +
         \sum_{\substack{ Y < p < R \\ p|u}}
            \frac{ \tau_{k-1}(u/p) }{\phi(a_ju/p)}  \pi^\flat(a_j N/p) 
            \label{E:Lemma5} \\
       &  + O(\tau_{k} (u))  + O\left( \tau_{k-1}(u) \Delta^*_\beta(AN;a_j u) \right) \notag \\
       &  + O\left( \tau_{k-1}(u) \sum_{p|u} \Delta^*(AN/p; a_ju/p) \right). \notag
\end{align}  
\end{lemma}

\begin{proof}
Assume that  $u|P_{\cL}(n)$ and let $u_i=(P_{\cL}(n),u)$.
Then $u=u_1 \ldots u_k$, where each $u_i|P_i(n)$. 
Moreover, we claim that this decomposition is unique because $\cL$ satisfies Hypothesis A. 
To justify this, assume that the decomposition is not unique. 
Then there must be some prime $p$ such that $p|L_i(n)$ and $p|L_j(n)$
for distinct values of $i$ and $j$.
We conclude that  $p|(a_ib_j-a_jb_i)$; this, however, contradicts Hypothesis A.

Therefore
\begin{equation} \label{E:51}
\sum_{\substack{N< n \le 2N \\ u|P_\cL(n)}} \beta(L_j(n)) =
\sum_{u_1 \ldots u_k = u} \sum_{\substack{ N < n \le 2N \\ u_i|L_i(n) \\ i=1,\ldots, k}} 
  \beta(L_j(n)).
\end{equation}
Set $L_j(n)=m$. 
Then $a_j N + b_j < m \le 2a_j N+ b_j$
Moreover, $u_j|m$. Now when $\beta(m)\ne 0$, 
$m$ has exactly one prime divisor $p$ with $p\le N^{1/2}$,
and all prime divisors of $u$ are less than $R \le N^{\vartheta/2}\le N^{1/2}$.
Consequently, we may assume that either $u_j=1$ or $u_j=p$ for some prime $p< R$.
In the latter case, our definition of $\beta$ implies that we may also assume  $p>Y$.

From our definition of $m$, we  also have  
\begin{equation*}
m \equiv b_j \pmod {a_j} \text{ and } a_i m \equiv a_i b_j - a_j b_i \pmod {u_i }
\text{ for } i \ne j.
\end{equation*}
We use the Chinese Remainder Theorem to combine these into one congruence
\begin{equation*}
m\equiv  m_0 \pmod {a_j u/u_j}.
\end{equation*}
Observe that $m_0$ is relatively prime to $a_ju/u_j$ by Hypothesis A;
the condition $(u,A)=1$ implies that $u$ is coprime to $a_i, a_j,$ and $a_ib_j-a_jb_i$. 

Now we fix values of $u_1, \ldots, u_k$. The inner sum of \eqref{E:51} is
\begin{equation*}
\sum_{\substack{ a_j N < m \le 2a_j N \\ m \equiv m_0 \pmod{a_j u/u_j}  }}
    \beta(m)
  + O(1).
\end{equation*}
Summing the error term over all values of $u_1, \ldots, u_k$ gives the
first error term in \eqref{E:Lemma5}.

Next, we consider the effect of different values of $u_j$.  First, we assume that $u_j=1$.
Then 
\begin{align} 
\sum_{\substack {a_j N < m \le 2a_j N \\ m \equiv m_0 \pmod{a_j u} }}
 \beta(m)   \label{E:uj=1}
 & = 
 \frac{\pi_{\beta, a_j u} (a_j N)}{\phi(a_j u)} +  \Delta_\beta(a_j N; a_j u, m_0) \\
 & = \frac{1}{\phi(a_j u)} \sum_{\substack {Y < p \le N^{1/2}  \\ p\nmid u} }\pi^\flat (a_j N/p) 
       +    \Delta_\beta(a_j N; a_j u, m_0).
       \notag
\end{align}

Now, assume that $u_j=p$ for some prime $p$, $Y < p < R$.
Let  $\bar{p}$ be the inverse of $p\pmod {a_j u/p}$. Then 
\begin{align}
\sum_{\substack {a_j N < m \le 2a_j N \\ m \equiv m_0 \pmod {a_j u/u_j}}  }  \beta(m) 
& =
         \sum_{\substack{ \frac{a_j N}{p} < p_2 \le \frac{2a_j N}{p} \\ 
                                                   p_2 \equiv m_0 \bar{p} \pmod {a_j u/p} } } 
                         \beta(pp_2)  \label{E:uj=p} \\
 &=  \pi^\flat(a_jN/p;a_ju/p,m_0\bar{p})   \notag \\
 &= \frac{\pi^\flat(a_jN/p)}{\phi(a_j u/p)}+ \Delta(a_jN/p; a_ju/p, \bar{p} m_0).  \notag  
 \end{align}

We  now sum \eqref{E:uj=1} and \eqref{E:uj=p}  over all choices of $u_1, u_2, \ldots, u_k$ 
with $u_1 u_2 \ldots u_k=u$ to finish the proof of the lemma.
\end{proof}

Returning to the estimate of $\cS_{1,j}$, we inject Lemma \ref{L:Lemma5} into 
\eqref{eq:50}. The first two error terms contribute
\begin{equation*}
\ll \sum_{d,e <R} |\lambda_d \lambda_e| \tau_{k} ([d,e]) 
    \left\{ 1+ \Delta^*_\beta (AN,a_j [d,e]) \right\}.
\end{equation*}
Using \eqref{eq:lambdabound}, we find that this is
\begin{equation*}
 \ll (\log R)^{2k} \sum_{r<R^2} \mu^2(r) (3k+3)^{\omega(r)}  \left\{ 1+ \Delta^*_\beta (AN,a_j r) \right\}.
\end{equation*}
By Lemma \ref{L:Lemma4}, 
there is some constant $C$ such that if $R\le N^{\vartheta/2} (\log N)^{-C}$, then 
the above is 
$\ll N$. 

The contribution of the third error term requires a slightly more elaborate analysis.
After using \eqref{eq:lambdabound}, we find that this contribution is
\begin{align*}
\ll & (\log R)^{2k} \sum_{r< R^2} \mu^2(r) (3k)^{\omega(r)} 
   \sum_{p|r; p<R} \Delta^*(AN/p, a_j r/p) \\
\ll & (\log R)^{2k} \sum_{p<R} 
          \sum_{m<R^2/p} \mu^2(m) (3k)^{\omega(m)} \Delta^*(AN/p,a_j m).
\end{align*}
We use Lemma \ref{L:Lemma3} to bound the innermost sum. 
If $R \le N^{\vartheta/2} (\log N)^{-C}$ for some sufficiently large $C$, then
the above is 
\begin{equation*}
\sum_{p<R} \frac{N}{p (\log N/p)} \ll N.
\end{equation*}

We close this section by updating our progress on $\cS_{1,j}$.
So far, we have
\begin{equation*} 
\cS_{1,j} = \frac{1}{\phi(a_j)}  \sumprime_{Y < p \le  N^{1/2} } 
   \pi^\flat \left( \frac{a_j N}{p} \right) T_p + O(N),
\end{equation*} 
where we define
\begin{equation} \label{E:TpDef}
T_p = 
   \sum_{\substack{ d,e \\ p \nmid [d,e]}}
      \frac{  \lambda_d \lambda_e \tau_{k-1}([d,e])}{\phi([d,e])}
     +
   \sum_{\substack{ d,e \\ p| [d,e]}}
      \frac{  \lambda_d \lambda_e \tau_{k-1}([d,e]/p)}{\phi([d,e]/p)}.
\end{equation}
Note that in the sum defining $T_p$, we implicitly have the conditions $d<R$ and $e<R$
because we are assuming that $\lambda_d=0$ if $d\ge R$. 
Therefore, if $p\ge R$, the second sum in \eqref{E:TpDef} is empty,
and the condition that $p\nmid[d,e]$ is vacuous. 
In other words, if $R \le p < N^{1/2}$, then 
\begin{equation} \label{E:TpRa}
T_p=  \sum_{d,e}
      \frac{  \lambda_d \lambda_e \tau_{k-1}([d,e])}{\phi([d,e])}
\end{equation}
However, when $p < R$, the sum $T_p$ is more complicated, and 
we will analyze this case in more detail in the next section.

Before closing this section,  we use the prime number theorem to write
\begin{equation*}
\pi^\flat \left( \frac{a_j N}{p} \right) = \frac{a_j N}{p \log N}  \alpha(p)  
 + O\left( \frac{N}{p (\log N)^2} \right),
\end{equation*}
where
\begin{equation*} \label{eq:alphaDef}
 \alpha(p) =  \frac{\log N}{\log(N/p)}.
\end{equation*}
Note that by Hypothesis A, $a_j$ and $A$ have exactly the same prime divisors.
Consequently, $a_j/\phi(a_j)= A/\phi(A)$, and
\begin{equation} \label{eq:S1j-1}
\cS_{1,j} =
\frac{A}{\phi(A)}  \frac{N}{\log N} 
\sumprime_{Y < p \le N^{1/2}} \frac{ \alpha(p) }{p} T_p 
+ O\left(N+ \frac{N}{(\log N)^2} \sumprime_{Y < p \le N^{1/2}} \frac{T_p}{p} \right).
\end{equation} 
 
\section{Evaluation of $T_p$ } \label{sec:Tp}

Analogous to the function $f$ defined in \eqref{E:fDef}, we define
\begin{equation} \label{E:f*Def}
f^*(d) = \frac{\phi(d)}{\tau_{k-1}(d)}
\end{equation}
whenever $d$ is squarefree and relatively prime to $A$.
We use this to define 
\begin{equation} \label{E:TdeltaDef}
    T_\delta = \sumprime_{d,e} 
               \frac{\lambda_{d}\lambda_{e} }{f^*([d,e,\delta]/\delta)}.
\end{equation}
When $\delta=p$, \eqref{E:TdeltaDef} reduces to the earlier definition of $T_p$.
We will analyze the more general quantity $T_\delta$; this provides additional insight at the cost of 
little extra complication of detail. 

An expression similar to $T_\delta$ occurs in Selberg's $\Lambda^2 \Lambda^-$ sieve.
See, for example, the last displayed equation on page 85 of Selberg\cite{Sel} or 
equation (1.9) on page 287 of Greaves\cite{Gr}. In our notation, those results 
can be stated as
\begin{equation*}
\sum_{d,e} \frac{\lambda_d \lambda_e}{f([d,e,\delta]/\delta)} =
\sum_{\substack{r \\ (r,\delta)=1}} 
     \frac{\mu^2(r)}{f_1(r)} 
      \left( \sum_{s|\delta} \mu(s) y_{rs} \right)^2.
\end{equation*}
Our next lemma is an analogue of this result with $f$ replaced by $f^*$.

\begin{lemma} \label{L:Tdelta}
If $\delta$ is squarefree and relatively prime to $A$, then
\begin{equation*}
T_\delta = 
\sumprime_{\substack{r \\ (r,\delta)=1}} \frac{\mu^2(r)}{f_1^*(r)} 
 \left( \sum_{s|\delta} \mu(s) y_{rs}^* \right)^2.
\end{equation*}
where 
\begin{equation*}
      f_1^*(d) =  \mu*f^*(d) = \prod_{p|d} \frac{p-k}{k-1}
\end{equation*}
whenever $d$ is squarefree and $(d,A)=1$,  and 
\begin{equation} \label{eq:yr*Def}
y_r^* = \frac{\mu^2(r)r}{\phi(r)} \sumprime_m \frac{y_{mr}}{\phi(m)}.
\end{equation}
\end{lemma}

\begin{proof}
Define $g^*_\delta(d)=g^*(d)$ by the relation
\begin{equation*}
    g^*(d)= f^*\left( \frac{d}{(d,\delta)} \right).
\end{equation*}
If $p$ is a prime, $p\nmid A$, then  
\begin{equation*}
    g^*(p)= 
    \begin{cases}
	 f^*(p) & \text{ if $p\nmid \delta $} \\
	 1 & \text{ if $p|\delta$. }
     \end{cases}
\end{equation*}
With this notation, we may write
\begin{equation*}
    T_\delta= \sumprime_{d,e} 
          \frac{\lambda_d \lambda_e g^*((d,e))}
	           {g^*(d) g^*(e)}
       = \sumprime_{d,e} \frac{\lambda_d \lambda_e}{g^*(d)g^*(e)}
              \sum_{\substack{r|d\\ r|e}} g^*_1(r),
\end{equation*}
where $g_1^*=g^* *\mu$.  
Note that 
\begin{equation*}
    g_1^*(p)= 
    \begin{cases} 
	 f^*_1(p) & \text{ if $p\nmid \delta$},\\
	 0 & \text{ if $p|\delta$.}
    \end{cases}
\end{equation*}
After changing the order of summation in the last sum, we find 
that 
\begin{equation*}
    T_\delta=\sumprime_{\substack{r\\ (r,\delta)=1}} g_1^*(r) 
      \left( \sumprime_{\substack{d \\ r|d }} \frac{\lambda_d}{g^*(d)} 
      \right)^2.
\end{equation*}
The condition that $(r,\delta)=1$ may be inserted because 
$g_1^*(r)=0$ if $(r,\delta)\ne 1$.

Define
\begin{equation*} 
    w_r^* = \mu(r) g_1^*(r) \sumprime_d \frac{\lambda_{dr}}{g^*(dr)}.
\end{equation*}
Then
\begin{equation}  \label{eq:Tdelta1}
    T_\delta=\sumprime_{\substack{r\\ (r,\delta)=1}} 
        \frac{\mu^2(r)}{g_1^*(r)} (w_r^*)^2.
\end{equation}

Assume henceforth that $(r,\delta)=1$.
Then
\begin{align*}
w_r^* = & \frac{\mu(r)g_1^*(r)}{g^*(r)} 
               \sumprime_d \frac{\lambda_{dr}}{g^*(d)} \\
      = & \frac{\mu(r)g_1^*(r)}{g^*(r)} 
               \sumprime_d \frac{\mu(dr) f(dr)}{g^*(d)}
                 \sumprime_t \frac{y_{drt}}{f_1(drt)} \\
      = & \frac{\mu^2(r) g_1^*(r) f(r)}{g^*(r) f_1(r)}
             \sumprime_m \frac{y_{mr}}{f_1(m)} 
              \sumprime_{d|m} \frac{\mu(d)f(d)}{g^* (d)}.
\end{align*}

We note that $g_1^*(r)=f_1^*(r)$ and 
$g^*(r)=f^*(r)$ because of our hypothesis $(r,\delta)=1$.
Thus
\begin{equation*}
\frac{g_1^*(r)f(r)}{g^*(r)f_1(r)} = 
\frac{f_1^*(r)f(r)}{f^*(r)f_1(r)} =
\frac{r}{\phi(r)}.
\end{equation*}

Next, we consider the sum
\begin{equation*}
\sumprime_{d|m} \frac{\mu(d)f(d)}{g^*(d)} =
\sumprime_{d|m} \frac{\mu(d)f(d)}
     {f^*\left( {\ds \frac{d}{(d,\delta)} }\right)}.
\end{equation*}
We write $d=d_1d_2$, with $(d_1,\delta)=1$ and $d_2|\delta$. 
The above sum is then 
\begin{equation} \label{E:Stuff1}
\sumprime_{\substack{d_1|m \\ (d_1,\delta)=1}}
 \frac{\mu(d_1) f(d_1)}{f^*(d_1)}
\sumprime_{\substack{d_2|m\\d_2|\delta}} \mu(d_2) f(d_2).
\end{equation}
The first factor in \eqref{E:Stuff1} is
\begin{equation*}
    \sumprime_{\substack{d_1|m \\ (d_1,\delta)=1}}
     \frac{\mu(d_1) f(d_1)}{f^*(d_1)}
  = \prod_{\substack{p|m\\p\nmid \delta}}
      \left( 1 - \frac{p(k-1)}{(p-1)k}\right)
= \frac{f_1(m/(m,\delta))}{\phi(m/(m,\delta))}.
\end{equation*}
The second factor in \eqref{E:Stuff1} is 
\begin{equation*}
\sum_{\substack{d_2|(m,\delta)} }
 \mu(d_2) f(d_2)=
\prod_{p|(m,\delta)} (1-f(p))=
\mu((m,\delta)) f_1((m,\delta)).
\end{equation*}
We conclude that the expression in \eqref{E:Stuff1} is
\begin{equation*}
\frac{f_1(m)}{\phi(m)} \mu((m,\delta)) \phi((m,\delta)).
\end{equation*}

Now 
\begin{equation*}
\mu((m,\delta)) \phi((m,\delta)) = 
\sum_{\substack{s|m\\s|\delta}} \mu(s) s,
\end{equation*}
so
\begin{equation*}
w_r^* = \frac{\mu^2(r) r}{\phi(r)} 
           \sum_m \frac{y_{mr}}{\phi(m)} 
             \sum_{\substack{s|m \\ s|\delta}} \mu(s) s.
\end{equation*}

The definition of $w_r^*$ depends on $r$ as well as $\delta$. 
Using the definition of $y_r^*$ given in \eqref{eq:yr*Def}, we find that 
\begin{equation*}
w_r^* = \sum_{s|\delta} \mu(s) y_{rs}^*.
\end{equation*}
Inserting this into \eqref{eq:Tdelta1} completes the proof of the lemma.
\end{proof}

When $\delta=p$, Lemma \ref{L:Tdelta} becomes
\begin{equation} \label{E:TpG}
T_p =  \sumprime_{\substack{r\\(r,p)=1}} 
          \frac{\mu^2(r)}{f_1^*(r)} (y_r^*-y_{rp}^*)^2.
\end{equation}
Now $y_{rp}^*=0$ if $r\ge R/p$, so 
\begin{equation*} \label{E:Tp1}
T_p =  \sumprime_{\substack{r < R/p \\(r,p)=1}} 
          \frac{\mu^2(r)}{f_1^*(r)} (y_r^*-y_{rp}^*)^2
         +  \sumprime_{\substack{R/p \le r < R \\(r,p)=1}} 
          \frac{\mu^2(r)}{f_1^*(r)} (y_r^*)^2.
\end{equation*}
When $p\ge R$, the second sum above is empty, and the condition
$(p,r)=1$ in the first sum is vacuous. 
In other words, if $p\ge R$, then
\begin{equation*} \label{E:TpRb}
T_p =  \sumprime_{r< R}
          \frac{\mu^2(r)}{f_1^*(r)} (y_r^*)^2 = T_1.
\end{equation*}
This is equivalent to  the observation that we made earlier in
\eqref{E:TpRa}.

Now we turn our attention to the sum 
\begin{equation*} 
\sumprime_{Y < p \le N^{1/2}} \frac{ \alpha(p)}{p} T_p
\end{equation*}
that appears in the main term of \eqref{eq:S1j-1}.
Using the above observations on $T_p$, we find that
\begin{equation}  \label{eq:MTDecomp}
\sumprime_{Y < p \le N^{1/2}} \frac{ \alpha(p)} {p} T_p = S_1+S_2+S_3,
\end{equation}
where 
\begin{align}
S_1= & \sumprime_{Y < p< R} \frac{ \alpha(p)}{p}
 \sumprime_{\substack{r< R/p \\ (r,p)=1}} 
    \frac{\mu^2(r)}{f_1^*(r)} (y_r^*-y^*_{rp})^2 ,
     \label{eq:S1Def} \\
S_2= & \sumprime_{Y < p< R} \frac{ \alpha(p)}{p}
 \sumprime_{\substack{R/p \le r < R \\ (r,p)=1}} 
     \frac{\mu^2(r)}{f_1^*(r)} (y_r^*)^2 ,
     \label{eq:S2Def} \\ 
S_3= & \sumprime_{R \le p< N^{1/2}} \frac{ \alpha(p)}{p}
      \sumprime_{r< R} 
         \frac{\mu^2(r)}{f_1^*(r)} (y_r^*)^2.
         \label{eq:S3Def}
\end{align}

\begin{lemma} \label{L:yr*Est}
 Assume that $r < R$, $(r,A)=1$, and $r$ is squarefree.
Let $y_r^*$ be as defined in \eqref{eq:yr*Def}. 
Then
\begin{equation*}
y_r^*= \frac{\phi(A)}{A} \fS(\cL) (\log R) \Pt\left( \frac{\log R/r}{\log R} \right)  + O(L(r)),
\end{equation*}
where
\begin{equation*}
L(r) = 1+ \sum_{p|r} \frac{\log p}{p}.
\end{equation*}
\end{lemma}

\begin{proof}
From \eqref{eq:yr*Def} and \eqref{E:yrChoice},  we see that
\begin{equation} \label{E:yr*Formula}
y_r^* = 
  \mu^2(r) \frac{r}{\phi(r)} \fS(\cL)
       \sum_{\substack{m \le R/r\\ (m,rA)=1} } \frac{\mu^2(m)}{\phi(m)}
        P\left( \frac{\log R/rm}{\log R} \right).
\end{equation}
We apply Lemma \ref{L:WirsingT} with 
\begin{equation*}
\gamma(p) = 
 \begin{cases} 
      1 & \text{ if $p\nmid rA$}, \\ 
      0 & \text{ if $p|rA$}.
  \end{cases}
\end{equation*}
Then $c_\gamma= \phi(rA)/rA$ and 
condition \eqref{E:Omega2} is satisfied with $\kappa=1$ and
\begin{equation*}
L = \sum_{p|rA} \frac{\log p}{p} + O(1).
\end{equation*}
We are regarding $A$ as fixed, so $L\ll L(r)$. 
Using Lemma \ref{L:WirsingT} with
\begin{equation*}
F(x) = P \left(x \frac{\log R/r}{\log R} \right),
\end{equation*}
we obtain
\begin{align*}
\sum_{\substack{ m < R/r \\ (m,rA)=1}} \frac{\mu^2(m)}{\phi(m)}  
  P\left( \frac{\log R/rm}{\log R} \right) =  
  \frac{\phi(rA)}{rA} (\log (R/r)) &
     \int_0^1 P\left( \frac{\log R/r}{\log R} (1-x) \right) dx  \\
    & + 
     O\left( \frac{\phi(r)}{r} L(r) \right).
\end{align*}

The desired results follows by making an appropriate change of variables in the integral on the 
right-hand side.
\end{proof}

\begin{lemma} \label{L:G*}
For $u\ge 1$, define
\begin{equation*}
G^*(u):= \sumprime_{r < u } \frac{\mu^2(r)}{f_1^*(r)}.
\end{equation*}
Then 
\begin{equation} \label{eq:G*est}
G^*(u) = \frac{A}{\phi(A)} \frac{(\log u)^{k-1}}{\fS(\cL) (k-1)!} + E^*(u),
\end{equation} 
where $E^*(u) \ll (\log (2u))^{k-2}$.
\end{lemma}

\begin{proof}
We apply  Lemma \ref{L:Wirsing} with
\begin{equation} \label{eq:L8gamma}
\gamma(p) = 
   \begin{cases}
   {\displaystyle  \frac{p(k-1)}{p-1} }& \text{ if $p\nmid A$,} \\ 0 & \text{if $p|A$,}\end{cases}
\end{equation}
and $\kappa=k-1$. 
As noted in \eqref{eq:2.2a}, every prime $p\le k$ divides  $A$, so
\begin{equation*}
\frac{\gamma(p)}{p} \le 1 - \frac{1}{k}
\end{equation*} 
for $p\nmid A$. Therefore, \eqref{E:Omega1} is satisfied with $A_1=k$. 
We are treating $A$ as fixed, so \eqref{E:Omega2} is satisfied with $L\ll 1$.
Moreover,
\begin{equation*}
c_\gamma= \prod_{p|A} \left( 1 - \frac{1}{p} \right)^{k-1} 
                    \prod_{p\nmid A} \left( 1-\frac{k-1}{p-1} \right)^{-1} \left( 1-\frac{1}{p} \right)^{k-1}
                  = \frac{A}{\phi(A)} \frac{1}{\fS(\cL)},
\end{equation*}
and \eqref{eq:G*est} follows  from Lemma \ref{L:Wirsing}.
\end{proof}

\begin{lemma} \label{L:Lemma8}
If  $p$ is prime, then  $T_p\ll (\log R)^{k+1}$.
\end{lemma}

\begin{proof}
From \eqref{E:TpG}, 
\begin{equation*}
T_p =  \sumprime_{\substack{r\\(r,p)=1}} 
          \frac{\mu^2(r)}{f_1^*(r)} (y_r^*-y_{rp}^*)^2.
\end{equation*}
We are regarding $P$ as fixed, so  Lemma \ref{L:yr*Est} implies that $y_j^*\ll \log R$ 
for any $j< R$, and $y_j^*=0$ if $j\ge R$.
Therefore
\begin{equation*}
T_p \ll (\log R)^2 \sumprime_{r < R} \frac{\mu^2(r)}{f_1^*(r)},
\end{equation*}
and the lemma now follows by using \eqref{eq:G*est} with $u=R$.
\end{proof}

Using the above lemma, we see that the second error term in \eqref{eq:S1j-1} is 
\begin{equation} \label{eq:FirstError}
\frac{N}{(\log N)^2} (\log R)^{k+1}  \sum_{p \le N^{1/2}} \frac{1}{p} 
\ll 
N (\log \log N) (\log N)^{k-1} .
\end{equation}

Combining  \eqref{eq:S1j-1}, \eqref{eq:FirstError}, and \eqref{eq:MTDecomp}, we now have
\begin{equation} \label{eq:S1j-2}
\cS_{1,j} =\frac{A}{\phi(A)}  \frac{N}{\log N} (S_1+S_2+S_3)
+ O\left( N (\log\log N) (\log N)^{k-1} \right). 
\end{equation}

To finish the proof of Theorem  \ref{T:S1Est}, 
we will show that when $i=1,2,$ or $3$,
\begin{equation*}
S_i =  \frac{\phi(A)}{A} \frac{\fS(\cL) (\log R)^{k+1} }{(k-2)!} J_i + O\left((\log \log R) (\log R)^k \right),
\end{equation*}
where $J_1, J_2, J_3$ are as defined in the statement of Theorem \ref{T:S1Est}.

\section{Completion of proof of Theorem \ref{T:S1Est}} \label{sec:7}

\begin{lemma} \label{L:S1Est}
Let $S_1$ be as defined in \eqref{eq:S1Def},
and let  $J_1$ be as defined in the statement of Theorem \ref{T:S1Est}.
Then 
\begin{equation*}
S_1 = \frac{\phi(A)}{A} \frac{\fS(\cL) (\log R)^{k+1}}{(k-2)!} J_1 + O\left( (\log \log R) (\log R)^k  \right).
\end{equation*}
\end{lemma}
 
\begin{proof} Assume that $r<R $, $r$ is squarefree, $p$ is a prime
with $p< R/r, (p,r)=1$, and $(pr,A)=1$.  By Lemma \ref{L:yr*Est},
\begin{equation*}
y_r^* - y_{rp}^* = 
    \frac{\phi(A)}{A} \fS(\cL) (\log R) 
      \int_{1- \frac{\log pr}{\log R}}^{1-\frac{\log r}{\log R}} P(x) dx + 
       O(L(r))
\end{equation*}
In  the above, we have used the simple observation that
\begin{equation*}
 L(rp)= L(r) + \frac{\log p}{p} \ll L(r)+1 \ll L(r).
\end{equation*}
Note also that
\begin{equation*} 
 (\log R) 
      \int_{1- \frac{\log pr}{\log R}}^{1-\frac{\log r}{\log R}} P(x) dx
      \ll \log p ,
\end{equation*}
and 
\begin{align}   
L(r) &  \le 1 + \sum_{p\le \log R} \frac{\log p}{p} + 
               \sum_{\substack{p|r \\ p > \log R}} \frac{\log p}{p}\label{eq:LrBound}\\
        &   \ll 1 + \log \log R + \frac{\log\log R}{\log R} \frac{\log R}{\log \log R}  
             \notag \\
          &  \ll \log \log R. \notag
\end{align}
In particular, $L(r) \ll \log p$ when $p > Y$. 
Therefore
\begin{align*}
 (y_r^* - y_{rp}^*)^2 = \frac{\phi(A)^2}{A^2} \fS(\cL)^2 (\log R)^2 & 
             \left(   \Pt\left( \frac{\log R/r}{\log R}\right)   -
                           \Pt\left( \frac{\log R/rp}{\log R}\right)  \right)^2 \\
                    &+ O\left (( \log  p) L(r) \right).
\end{align*}
We use this in the definition of $S_1$ to obtain  
\begin{align*}
S_1 = & 
     \frac{\phi(A)^2}{A^2}  \fS(\cL)^2 (\log R)^2  \\
     & \text{\hskip .25in}  
  \sumprime_{Y < p<R} \frac{ \alpha(p)}{p} 
    \sumprime_{\substack{r<R/p \\ (r,p)=1}} \frac{\mu^2(r)}{f_1^*(r)}  
     \left( \Pt\left( \frac{\log R/r}{\log R}\right)- 
                      \Pt\left( \frac{\log R/rp}{\log R}\right)
                 \right)^2
       \\
     & \text{\hskip  1in}
        + O\left( \sum_{p<R}  \frac{1}{p} 
                            \sumprime_{r<R/p}  \frac{\mu^2(r)}{f_1^*(r)} (\log p) L(r) 
                                \right)    \\
          = & S_{11} + O(S_{12}),
\end{align*}
say.

For $S_{12}$, we   reverse the order of summation and use \eqref{eq:G*est} to obtain 
\begin{align}
S_{12} &  = \sumprime_{r< R} \frac{\mu^2(r)}{f_1^*(r)} L(r) \sum_{p< R/r} \frac{\log p}{p}     
               \ll (\log\log R) (\log R) \sumprime_{r<R} \frac{\mu^2(r)}{f_1^*(r)}  \label{eq:S12} \\
                & \ll (\log \log R) (\log R)^k. \notag
\end{align}

Now we consider $S_{11}$.
We write this as $ S_{13}- S_{14}$, where
\begin{align} \label{eq:S13}
S_{13} =    \frac{\phi(A)^2}{A^2} &  \fS(\cL)^2 (\log R)^2  \\
  &  \sumprime_{Y< p<R} \frac{\alpha(p)}{p} 
    \sumprime_{r<R/p} \frac{\mu^2(r)}{f_1^*(r)}   
      \left( \Pt\left( \frac{\log R/r}{\log R}\right)- 
                      \Pt\left( \frac{\log R/rp}{\log R}\right)
                 \right)^2, \notag
\end{align}   
and $S_{14}$ is the same sum with the extra condition that $p|r$.

For $S_{14}$, we note that
\begin{equation*}
\Pt\left( \frac{\log R/r}{\log R} \right) - \Pt\left(\frac{\log R/rp}{\log R} \right)
= \int_{\frac{\log R/rp}{\log R}}^{\frac{\log R/r}{\log R}} P(t) dt \ll
\frac{\log p}{\log R}.
\end{equation*}
We also note that $f_1^*(p) = (p-k)/(k-1) \gg p$.
Making the change of variables $r=mp$, we get
\begin{align*}  \label{eq:S14}
S_{14} \ll &(\log R)^2 \sumprime_{Y<p<R} \frac{1}{p^2} 
                 \sumprime_{m<R/p} \frac{\mu^2(m)}{f_1^*(m)} 
                  \left( \frac{\log p}{\log R} \right)^2 \\
             \ll & (\log R)^{k-1} \sum_{p<R} \frac{(\log p)^2}{p^2}
\end{align*}
by Lemma 8. The last sum converges, so
\begin{equation} \label{eq:S14}
S_{14} \ll (\log R)^{k-1}.
\end{equation}

For $S_{13}$, we evaluate the inner sum using Lemma \ref{L:WirsingT}
with $z=R$, $g(d)=1/f_1^*(d)$, $\kappa=k-1$, $\gamma$ as defined in \eqref{eq:L8gamma},
and
\begin{equation*}
 F\left( \frac{\log R/r}{\log R} \right)
  = 
   \begin{cases}
   \left\{ \Pt\left( \frac{\log R/r}{\log R} \right) - \Pt\left(\frac{\log R/rp}{\log R} \right)\right\}^2 &
       \text{ if $r< R/p$,}\\
     0 & \text{ if $R/p \le r < R$.}
  \end{cases}
\end{equation*}
If we set $y=\log p/\log R$ and $x=\log r/\log R$, then the last is equivalent to 
\begin{equation*}
F(1-x)= F_p(1-x) = 
 \begin{cases}
  \left( \Pt(1-x) - \Pt(1-x-y) \right)^2 & \text{ if $x < 1-y$},\\
  0 & \text{ if $1-y\le x < 1$}.
 \end{cases}
\end{equation*}
Making the substitution $w=1-x$, we see that this is the same as 
\begin{equation*}
F(w)= F_p(w) = 
 \begin{cases}
  \left( \Pt(w) - \Pt(w-y) \right)^2 & \text{ if $y \le w \le 1$},\\
  0 & \text{ if $0 \le w < y$}.
 \end{cases}
\end{equation*}
From Lemma \ref{L:WirsingT}, we find that 
\begin{align}
 \sumprime_{r<R/p} \frac{\mu^2(r)}{f_1^*(r)}   &
      \left( \Pt\left( \frac{\log R/r}{\log R}\right)- 
                      \Pt\left( \frac{\log R/rp}{\log R}\right)
                 \right)^2 = \label{eq:S13Inner} \\
   &   \frac{A}{\phi(A)}
     \frac{(\log R)^{k-1}}{\fS(\cL) (k-2)!} V_1\left(\frac{\log p}{\log R}\right)+
     O\left( M(F_p) (\log R)^{k-2} \right). \notag
\end{align}
where
\begin{equation*}
V_1(y) = \int_0^{1-y} \left\{ \Pt(1-x) - \Pt(1-x-y) \right\}^2 x^{k-2} dx.
\end{equation*}
Observe that if $y\le x \le 1$, then
\begin{equation*}
|F_p(x)| = \left(  \int_{x-y}^{x} P(t) dt \right)^2 \le y^2 \sup_{t\in [0,1]} |P(t)| \ll 1,
\end{equation*}
where the implied constant depends on  $P$ but not on $p$. 
Similarly, $F'_p(x) \ll 1$, and therefore $M(F_p) \ll 1$ uniformly in $p$. 
The error term in \eqref{eq:S13Inner} thus contributes
\begin{equation} \label{eq:S13ET}
\ll (\log R)^k \sum_{p<R} \frac{1}{p} \ll (\log \log R)(\log R)^k
\end{equation}
to $S_{13}$. 
Incorporating the contribution of the main term from \eqref{eq:S13Inner},
we now have
\begin{equation} \label{eq:S13a1}
 S_{13}=
  \frac{\phi(A)}{A}
     \frac{\fS(\cL) (\log R)^{k+1}}{(k-2)!} 
     \sumprime_{Y <p < R}  \frac{\alpha(p)}{p} V_1\left(\frac{\log p}{\log R}\right)
     + O\left( (\log \log R) (\log R)^k \right).
\end{equation}

Now let  $Z(u)$ be defined by the relation
\begin{equation} \label{eq:ZDef}
\sum_{p\le u } \log p = u + Z(u).
\end{equation}
From the classical form of the prime number theorem, we know that 
\begin{equation*}  
Z(u) \ll u \exp(-c\sqrt{\log u})
\end{equation*}
for some absolute constant $c$. 
Therefore the sum in \eqref{eq:S13a1} is 
\begin{equation} \label{eq:S13a2}
 \int_Y^R   \alpha(u)  V_1\left(\frac{\log u}{\log R}  \right) \frac{du}{u\log u}   
   + \int_Y^R   \alpha(u) V_1\left(\frac{\log u}{\log R} \right) \frac{ dZ(u)} {u\log u} .
\end{equation}
In the first integral, we make the change of variable $u=R^y$, and we set
\begin{equation} \label{E:bDef}
b= \frac{\log N}{\log R}
\end{equation}
to obtain
\begin{equation} \label{eq:S13a21}
\int_Y^R   \alpha(u)  V_1\left(\frac{\log u}{\log R}  \right) \frac{du}{u\log u}  
=   \int_{b\eta}^1 \frac{b}{y(b-y)} V_1(y) dy =J_1.
 \end{equation}
 Note that we have used the fact that $\log Y/\log N=\eta$.
 Comparing the definitions of $b$ and $B$ (see \eqref{eq:1.18}), we see 
 that 
 \begin{equation*}
 b= \left(  \frac{\vartheta}{2} - \frac{C\log \log N}{\log N} \right)^{-1}  
   = B+O(\log \log R/\log R).
 \end{equation*}
We may therefore replace $b$ by $B$ on the right-hand side of 
 \eqref{eq:S13a21} at the cost of an error term $O(\log \log R/\log R)$. 
The first integral in 
 \eqref{eq:S13a2} is thus 
 \begin{equation} \label{eq:S13a3}
=  \int_{B\eta}^1 \frac{B}{B-y} V_1(y) dy + O\left( \frac{\log \log R}{\log R} \right) 
 =
 J_1 + O\left( \frac{\log \log R}{\log R} \right).
 \end{equation}

We write the second integral in \eqref{eq:S13a2} as 
\begin{equation*} 
\int_Y^R F_1(u) dZ(u) = F_1(R)Z(R)-F_1(Y)Z(Y) - \int_Y^R Z(u)F_1'(u) du,
\end{equation*}
where
\begin{equation*}
F_1(u) = \frac{\alpha(u)}{u\log u} V_1\left( \frac{\log u}{\log R} \right).
\end{equation*}
Now $V_1(y)\ll y$, so 
\begin{equation*}
F_1(u) \ll \frac{1}{u \log R}.
\end{equation*}
Moreover,
\begin{align*}
F_1'(u) = & F_1(u) \frac{d}{du} \log F_1(u) \\
        =  & F_1(u) 
                 \left\{ \frac{1}{u\log (N/u)} -
                             \frac{V_1'}{V_1} \left( \frac{\log u}{\log R} \right) \frac{1}{u\log R} -
                             \frac{1}{u} -\frac{1}{u\log u}
                  \right\},
\end{align*}
so 
\begin{equation*}
|F_1'(u)| \ll \frac{|F_1(u)|}{u} \ll \frac{1}{u^2 \log R}.
\end{equation*}    

Therefore
\begin{equation*} 
\int_Y^R Z(u) F_1'(u) du \ll   \int_Y^R \frac{\exp(-c\sqrt{\log u})}{u\log R} du 
                                       \ll  (\log R)^{-1}
\end{equation*} 
We also note that 
\begin{equation*}
|F_1(R)Z(R)|+|F_1(Y)Z(Y)|  \ll (\log R)^{-1}.
\end{equation*}
 
From the last two estimates, \eqref{eq:S13a1}, and \eqref{eq:S13a3}, 
we conclude that 
\begin{equation*}
S_{13} = \frac{\phi(A)}{A}  \frac{\fS(\cL) (\log R)^{k+1} } {(k-2)!} J_1 +
O\left( (\log \log R) (\log R)^k \right).
\end{equation*}
We combine this with  
\eqref{eq:S12} and \eqref{eq:S14}
to complete the proof.
\end{proof}

\begin{lemma} \label{L:S2Est}
Let $S_2$ be as defined in \eqref{eq:S2Def},
and let  $J_2$ be as defined in the statement of Theorem \ref{T:S1Est}.
Then 
\begin{equation*}
S_2 = \frac{\phi(A)}{A} \frac{\fS(\cL) (\log R)^{k+1}}{(k-2)!} J_2+ 
    O\left( (\log \log R)  (\log R)^k  \right).
\end{equation*}
\end{lemma}

\begin{proof} 
From Lemma \ref{L:yr*Est} and \eqref{eq:LrBound}, we see that
\begin{equation} \label{E:yr*2}
(y_r^*)^2 = \frac{\phi(A)^2}{A^2} \fS(\cL)^2 (\log R)^2 \Pt\left( \frac{\log R/r}{\log R} \right)^2 +
                    O(L(r) \log R).
\end{equation}

Therefore
\begin{align}
S_2 = & \frac{\phi(A)^2}{A^2}
              \fS(\cL)^2 (\log R)^2 \sumprime_{Y < p < R} \frac{\alpha(p)}{p} 
                 \sumprime_{\substack{ R/p \le r < R \\ (r,p)=1}} 
                     \frac{\mu^2(r)}{f_1^*(r)} \Pt\left( \frac{\log R/r}{\log R} \right)^2
                     \label{E:S2Decomp}  \\
            & \text{\quad} 
                   + O\left( (\log R) \sumprime_{Y < p < R} \frac{1}{p} 
                                  \sumprime_{R/p \le r < R } \frac{\mu^2(r)}{f_1^*(r)} L(r) \right)\notag\\
         = & S_{21} + O(S_{22}),\notag
\end{align}
say. 

We first consider $S_{22}$. From the above definition, we see that
\begin{equation*}
S_{22} \ll  (\log \log R) (\log R) \sumprime_{r<R}  \frac{\mu^2(r)}{f_1^*(r)} L(r).
\end{equation*}
Now 
\begin{equation*} 
 \sumprime_{r<R}  \frac{\mu^2(r)}{f_1^*(r)} L(r) =                          
            \sumprime_{r< R} \frac{\mu^2(r)}{f_1^*(r)}
                     + \sumprime_{r< R} \frac{\mu^2(r)}{f_1^*(r)}\sum_{p|r}  \frac{\log p}{p},
\end{equation*}
and
\begin{align*}
 \sumprime_{r< R} \frac{\mu^2(r)}{f_1^*(r)}  \sum_{p|r} \frac{\log p}{p} 
  = & \sumprime_{p< R} \frac{\log p}{pf_1^*(p)} 
              \sumprime_{t<R/p; (t,p)=1} \frac{\mu^2(t)}{f_1^*(t)} \\
  \ll & \sum_{p<R} \frac{\log p}{p^2} (\log R)^{k-1} \ll (\log R)^{k-1}.
 \end{align*}
 Therefore
 \begin{equation} \label{E:L(r)Sum}
  \sumprime_{r<R}  \frac{\mu^2(r)}{f_1^*(r)} L(r) \ll (\log R)^{k-1},
\end{equation}
and 
 \begin{equation}\label{E:S22}
 S_{22} \ll    (\log \log R) (\log R)^{k}.
 \end{equation}

Now $S_{21}=S_{23}-S_{24}$, where
\begin{equation} \label{E:S23Def}
S_{23} = \frac{\phi(A)^2}{A^2} \fS(\cL)^2 (\log R)^2 
                  \sumprime_{Y < p < R} \frac{\alpha(p)}{p}
                        \sumprime_{R/p \le r < R} \frac{\mu^2(r)}{f_1^*(r)} 
                            \Pt\left(\frac{\log R/r}{\log R}\right)^2,
\end{equation}
and $S_{24}$ is the same sum with the extra condition that $p|r$. 

For $S_{24}$, 
we begin by noting that $\Pt(y) \ll y$. Therefore, if $R/p \le r <R$, then
\begin{equation*}
\Pt\left( \frac{\log R/r}{\log R} \right)  \ll \left( \frac{\log p}{\log R} \right)^2.
\end{equation*}
Consequently,
\begin{equation} \label{E:S24}
S_{24} \ll (\log R)^2 \sum_{p<R} \frac{(\log p)^2}{p^2} 
                  \sumprime_{t<R} \frac{\mu^2(t)}{f_1^*(t)} 
            \ll (\log R)^{k-1}.
\end{equation}
                          
Using Lemma \ref{L:WirsingT}, we find that the innermost sum in $S_{23}$ is 
\begin{equation*}
  \frac{A}{\phi(A)}
     \frac{(\log R)^{k-1}}{\fS(\cL) (k-2)!} V_2\left(\frac{\log p}{\log R}\right)+
     O\left( (\log R)^{k-2} \right),
\end{equation*}
where
\begin{equation*}
V_2(y) = \int_{1-y}^1 \Pt(1-x)^2 x^{k-2} dx.
\end{equation*}
Inserting this into \eqref{E:S23Def}, we find that 
\begin{equation} \label{E:S23}
S_{23} = \frac{\phi(A)}{A} \frac{ \fS(\cL) (\log R)^{k+1} }{(k-2)!} 
       \sumprime_{Y < p < R} \frac{\alpha(p)}{p} V_2\left( \frac{\log p}{\log R}\right) +
       O\left( (\log\log R) (\log R)^k\right).
\end{equation}
 
The sum in the main term is 
\begin{equation} \label{eq:S23p}
\int_Y^R \alpha(u) V_2\left( \frac{\log u}{\log R}\right) \frac{du}{u\log u} +
\int_Y^R \alpha(u) V_2\left( \frac{\log u}{\log R} \right) \frac{dZ(u)}{u\log u},
\end{equation}       
where $Z(u)$ was defined in \eqref{eq:ZDef}. 
In the first integral, we let $u=R^y$ to obtain
\begin{equation*}
\int_Y^R \alpha(u) V_2\left( \frac{\log u}{\log R}\right) \frac{du}{u\log u}=
\int_{b\eta}^1 \frac{b}{y(b-y)} V_2(y) dy,
\end{equation*}
where $b=\log N/\log R$, as defined in \eqref{E:bDef}. 
As in the proof of Lemma \ref{L:S1Est}, we may
replace $b$ by $B$ at the cost of an error term 
$O(\log \log R/\log R)$; therefore,
\begin{equation*} 
\int_Y^R \alpha(u) V_2\left( \frac{\log u}{\log R}\right) \frac{du}{u\log u} =
J_2 + O\left( \frac{\log \log R}{\log R} \right).
\end{equation*}

The second integral in \eqref{eq:S23p} may be written as 
\begin{equation*}
\int_Y^R F_2(u) dZ(u)
\end{equation*}
where
\begin{equation*}
F_2(u) = \frac{\alpha(u)}{u\log u} V_2 \left( \frac{\log u}{\log R} \right).
\end{equation*}        
We estimate this by using the argument following \eqref{eq:S13a3}, but with 
$F_1$ and $V_1$ replaced by $F_2$ and $V_2$.    
Note that $V_2(y) \ll y$, so $F_2(u) \ll (u \log R)^{-1}$. 
The end result is that
\begin{equation*}
\int_Y^R F_2(u) dZ(u) \ll (\log R)^{-1}.
\end{equation*}
We combine the above estimates to get
\begin{equation*}
 \sumprime_{Y < p < R} \frac{\alpha(p)}{p} V_2\left( \frac{\log p}{\log R}\right) =
 J_2 + O\left(\frac{\log \log R} {\log R} \right).
 \end{equation*}
 The proof of the lemma is completed by combining this with
 \eqref{E:S2Decomp}, \eqref{E:S22}, \eqref{E:S24}, and \eqref{E:S23}.              
\end{proof}

\begin{lemma} \label{L:S3Est}
Let $S_3$ be as defined in \eqref{eq:S3Def},
and let  $J_3$ be as defined in the statement of Theorem \ref{T:S1Est}.
Then 
\begin{equation*}
S_3 = \frac{\phi(A)}{A} \frac{\fS(\cL) (\log R)^{k+1}}{(k-2)!} J_3 + 
 O\left( (\log \log R) (\log R)^k  \right).
\end{equation*}
\end{lemma}

\begin{proof} 
$S_3$ is a product of two sums. 
Using \eqref{E:yr*2}, we see that the second sum is
\begin{align*}
\sumprime_{r<R} \frac{\mu^2(r)}{f_1^*(r)} (y_r^*)^2 
= \frac{\phi(A)^2}{A^2} \fS(\cL)^2 (\log R)^2  &
     \sumprime_{r<R} \frac{\mu^2(r) }{f_1^*(r)} 
          \Pt\left( \frac{\log R/r}{\log R} \right)^2 \\
& + O\left( (\log R) \sumprime_{r<R} \frac{\mu^2(r)}{f_1^*(r)} L(r)\right).
 \end{align*}
We use Lemma \ref{L:WirsingT} for the main term and \eqref{E:L(r)Sum} for 
the error term.
Therefore
\begin{equation} \label{E:S3-2}
\sumprime_{r<R} \frac{\mu^2(r)}{f_1^*(r)} (y_r^*)^2 =
\frac{\phi(A)}{A} \frac{ \fS(\cL) (\log R)^{k+1}}{(k-2)!} \int_0^1 \Pt(1-x)^2 x^{k-2} dx + 
 O\left( (\log R)^k \right).
 \end{equation}
 
 The first sum in the definition of $S_3$ is 
 \begin{equation} \label{eq:S3p}
 \int_{R}^{N^{1/2}} \alpha(u) \frac{du}{u\log u} + \int_R^{N^{1/2}} \alpha(u) \frac{dZ(u)}{u\log u}.
 \end{equation}
 In the first integral, we set $u=R^y$ to get
 \begin{equation*}
 \int_{R}^{N^{1/2}} \alpha(u) \frac{du}{u\log u} = \int_1^{b/2} \frac{b}{y(b-y)} dy.
 \end{equation*}
 As in the proofs of Lemma \ref{L:S1Est} and Lemma \ref{L:S2Est}, 
 we may replace $b$ by $B$ at the cost of a small error term, and 
 therefore 
 \begin{equation*}
  \int_{R}^{N^{1/2}} \alpha(u) \frac{du}{u\log u}  = 
   \int_1^{B/2} \frac{B}{y(B-y)} dy + O \left( \frac{\log \log R}{\log R} \right).
\end{equation*}

Letting $F_3(u)=\alpha(u)/(u\log u)$, we see that the second integral 
in \eqref{eq:S3p} is 
 \begin{align*}
 \int_R^{N^{1/2}} & F_3(u)  dZ(u)   \\
        &  \ll |F_3(R) Z(R)| + |F_3(N^{1/2}) Z(N^{1/2}) | + 
                   \int_R^{N^{1/2}}  |F'_3(u)| \exp(-c\sqrt{\log u}) du \\
         & \ll (\log R)^{-1}.
\end{align*}
Therefore 
\begin{equation*}
\sumprime_{R \le p < N^{1/2} }\frac{ \alpha(p)}p =  \int_1^{B/2} \frac{B}{y(B-y)} dy + 
        O\left( (\log R)^{-1} \right).
\end{equation*}
We combine this with \eqref{E:S3-2} to complete the proof.
\end{proof}

Theorem 8 now follows by combining the previous three lemmas and \eqref{eq:S1j-2}.

We close this section by giving, as promised earlier,  a short sketch of the proof of Theorem \ref{th:9}. 
The left-hand side of the conclusion is
\begin{equation}  \label{E:9-1}
= \sum_{d,e} \lambda_d \lambda_e Q_j([d,e]),
\end{equation}
where
\begin{equation*}
Q_j(u) = \sum_{\substack{ N < n \le 2N \\ u|P_{\cL}(n) }} \varpi(L_j(n)).
\end{equation*}
This last sum can be evaluated in the same way as the related sum $M_j(u)$ considered 
in Lemma \ref{L:Lemma5}. The evaluation is simpler because only the case $u_j=1$ occurs in 
this instance. The final result is 
\begin{equation*}
Q_j(u) = \frac{\tau_{k-1}(u)}{\phi(a_j u)} \pi^\flat (a_j N) + O(\tau_{k-1}(u) \Delta^*(AN,u) ) .
\end{equation*}
We insert this into \eqref{E:9-1} and use the Bombieri-Vinogradov theorem to handle the error
terms. The main term is 
\begin{equation*}
\frac{\pi^\flat (a_j N)}{\phi(a_j)} T_1 = \frac{A}{\phi(A)} \frac{N}{\log N} T_1 + O(T_1 N ( \log N)^{-2} ), 
\end{equation*}
where $T_1$  is given by  \eqref{E:TdeltaDef} with $\delta=1$. By Lemma \ref{L:Tdelta}, $T_1$ is equal to the 
sum considered in \eqref{E:S3-2},  and the proof is completed by appealing to the formula there.

\section{Proofs of Theorems 1--3 and Corollaries}\label{sec:5}

For the proof of Theorem 1, we use \eqref{eq:cSDef}. For our choice of $P$, we take
$\ell=[\sqrt k]$, and 
\begin{equation*} 
P(x) = \frac{x^\ell}{\ell!} , {\quad} 
\Pt(x) = \frac{x^{\ell + 1}}{(\ell + 1)!}.
\end{equation*} 
We take $Y=1$ in the definition of $\beta$; therefore $\eta=0$.  

From Theorems \ref{T:S0Est} and \ref{T:S1Est}, we see that 
(cf. \eqref{eq:SDecomp})
\begin{equation*}
\cS \asymptotic
 \fS(\cL) N (\log R)^k J,
\end{equation*}
where
\begin{equation}  \label{eq:101a}
J = 
 \left\{ \frac{k}{B} \frac{(J_1+ J_2+ J_3)}{(k-2)!} -\nu \frac{J_0}{(k-1)!} \right\}.
 \end{equation}
Next, we write
 \begin{equation} \label{eq:101b}
 \frac{ J_1+ J_2}{(k-2)!} = J_4 + J_5 + J_6,
 \end{equation}
 where 
\begin{align}
J_4 &= 
   \int\limits^1_{0} 
      \left(\frac1{y} + \frac1{B -y}\right) 
    \int\limits^1_0 
        \widetilde P(1 - x)^2 \frac{x^{k -2}}{(k - 2)!}\, 
           dxdy, \nonumber \\
J_5 &= 
      -2 \int\limits^1_{0}
        \left(\frac1{y} + \frac1{B -y}\right) 
          \int\limits^{1 - y}_0 
            \widetilde P(1 - x) \widetilde P(1 - y - x) 
            \frac{x^{k - 2}}{(k - 2)!} \, dxdy,
\label{eq:102}\\
J_6 &= 
      \int\limits^1_{0}
         \left(\frac1{y} + \frac1{B -y}\right) 
        \int\limits^{1 - y}_0 \widetilde P(1 - x - y)^2
          \frac{x^{k - 2}}{(k - 2)!}\, dx dy.
\nonumber
\end{align}
In  $J_5$, we can write 
\begin{equation}
(1 - x)^{\ell + 1} = (1 - y - x)^{\ell + 1} + \sum^{\ell + 1}_{j= 1} 
   \binomial{\ell+1}{j}  
           y^j (1 - y - x)^{\ell + 1 - j}
\label{eq:103}
\end{equation}
and denote the corresponding integrals by $J^{(0)}_5$ and
$J^{(1)}_5$, resp.
The terms $J_4$, $J^{(0)}_5$ and $J_6$ will contribute to the
main term, $J^{(1)}_5$ to the secondary term.

We will often use the evaluation $(m, n \in \mathbb Z^+)$
\begin{equation}
\int\limits^1_0 x^m(1 - x)^n dx = \frac{m!n!}{(m + n + 1)!},
\label{eq:104}
\end{equation}
which is a special case of a standard formula for the Euler beta
function (see e.g.\ Karacuba \cite[p.~46]{Karac}).
For later convenience, we define%
\footnote{$A(k,\ell)$ should not be confused with the quantity $A$ defined in \eqref{eq:2.1}.}
\begin{equation*}
A(k,\ell) = \binomial{2\ell+2}{\ell+1} \frac{1}{(k+2\ell+1)!}.
\end{equation*}

Using \eqref{eq:104} we obtain
\begin{align}
J^{(0)}_5 + J_6
&= 
-J_6= -\int\limits^1_{0}  \left(\frac1{y} + \frac1{B - y}\right) 
\binomial{2\ell+2}{\ell+1} 
\frac{(1 - y)^{k + 2\ell +1}}{(k + 2\ell + 1)!} \, dy
\label{eq:105}\\
&= -A(k,\ell) \int\limits^1_{0} \left(\frac1{y} + \frac1{B - y}\right)(1 - y)^{k + 2\ell + 1} \, dy,
\nonumber
\end{align}
\begin{equation}
J_4 = 
A(k,\ell)  \int\limits^1_{0} \left(\frac1{y} + \frac1{B - y} \right) dy, 
\label{eq:106}
\end{equation}
\begin{align}
J_7  &= J_4 + J^{(0)}_5+ J_6 = 
A(k,\ell) \int\limits^1_{0}\left(\frac1{y} + \frac1{B - y} \right)
\left(1 - (1 - y)^{k + 2\ell + 1}\right)\, dy 
\label{eq:107}\\
&= A(k,\ell) \int\limits^1_{0} 
\Biggl(\frac1{y} y \sum^{k + 2\ell}_{j =0} (1 - y)^j + 
\frac1{B - y} - \frac{(1 - y)^{k + 2\ell + 1}}{B- y} \Biggr) dy.
\nonumber
\end{align}
With the notation
\begin{equation}
L(n) = \sum^n_{i = 1} \frac1{i} = \log n + \gamma + O\Bigl(\frac1{n}\Bigr),
\label{eq:108}
\end{equation}
we obtain by $B \geq 2$,
\begin{align}
J_7 &= A(k,\ell)  
\left(
   \sum^{k + 2\ell}_{j = 0} \int\limits^1_0 (1 -y)^j \, dy +
     \int\limits^1_0 \frac{dy}{B - y} - 
      \int\limits^1_0 \frac{(1 - y)^{k + 2\ell + 1}}{B - y} \, dy 
\right) 
\label{eq:109}\\
&= A(k,\ell)  
   \left( L(k + 2\ell + 1) + \log \frac{B}{B - 1} + O( 1/k) \right) 
\nonumber \\
& = A(k,\ell) 
 \left(\log k + \gamma +\log \frac{B}{B - 1} + O(1/\sqrt{k}) \right).
\nonumber
\end{align}

Since in the term $J^{(1)}_5$ the factor $y$ appears, we can
directly work with $J^{(1)}_5$, and we get 
\begin{align}
-\frac{J^{(1)}_5}{2} &= 
\int\limits^1_0 \sum^{\ell + 1}_{j =1} \frac{y^{j - 1} + 
 \frac{y^j}{B - y}}{(\ell + 1)! (\ell + 1 -j)! j!} 
  \int\limits^{1 - y}_0 \frac{(1 - y - x)^{2\ell + 2 - j} x^{k - 2}}{(k - 2)!} \, dx dy
\label{eq:110}\\
&=  \sum^{\ell + 1}_{j = 1} \int\limits^1_0
\frac{(y^{j - 1} + O(y^j))(1 - y)^{k + 2\ell + 1 - j}
 (2 \ell + 2 - j)!}{(\ell + 1)! (\ell + 1 - j)! j!(k + 2\ell + 1 - j)!} \,dy \nonumber\\
&= A(k,\ell) \sum^{\ell + 1}_{j = 1} \frac{(\ell + 1) \dots 
(\ell + 1 -(j - 1))}{(2\ell + 2) \dots (2\ell + 2 - (j - 1))} 
\left(\frac1j + O \left(\frac1{k + 2\ell + 2}\right)\right) \nonumber\\
&= A(k,\ell) \left(C(\ell) + O ( 1/k ) \right).
\nonumber
\end{align} 
By $\log (1 - x) = - \sum\limits^\infty_{j = 1} x^j / j$ we have 
\begin{equation}
C(\ell) = \sum^\infty_{j = 1} \frac1{2^j} \cdot \frac1j + o(1) =
-\log \left(1 - \frac12\right) + o(1) = \log 2 + o(1)
\label{eq:111}
\end{equation}
as  $\ell \to \infty$, so
\begin{equation*}
J_5^{(1)} = -\log 4 + o(1).
\end{equation*}

Finally,
\begin{equation}
\frac{J_3}{(k - 2)!} = \int\limits^{B/2}_1 \left(\frac1{y} +
\frac1{B - y} \right) \cdot \frac{(2\ell + 2)!}{((\ell + 1)!)^2} \frac{dy}{(k + 2\ell + 1)!} 
= A(k,\ell) \log(B - 1).
\label{eq:112}
\end{equation}

Summarizing \eqref{eq:101b}--\eqref{eq:112}, we obtain 
\begin{align}
\frac{J_1 + J_2 + J_3}{(k - 2)!} &=
 J_7 + J^{(1)}_5 +\frac{J_3}{(k - 2)!} + o(1)
\label{eq:113}\\
&= A(k,\ell)  \left(\log k + \gamma + \log B - \log 4 + o(1) \right)
\nonumber\\
&= A(k,\ell)  \left(\log \frac{B e^\gamma k}{4} + o(1)\right).
\nonumber
\end{align}

From  \eqref{eq:104}, we  deduce that
\begin{equation*}
\frac{J_0}{(k-1)!} = \binom{2\ell}{\ell} \frac{1}{(k+2\ell)!} = A_0(k,\ell),
\end{equation*}
say. 

Returning to  \eqref{eq:101a}, we find that
\begin{align*}
J & =
\frac{k}{B} A(k,\ell)  \log \left(\frac{B e^\gamma k}{4} \right) - \nu A_0(k,\ell)
 + o(kA(k,\ell))  \\
& =  A_0(k,\ell) 
           \left(
                      \frac{2 \left(2- \frac1{\ell + 1} \right)k}{B(k + 2\ell + 1)}
                         \log \left(\frac{Be^\gamma k}{4}\right) - \nu + o(1) \right) \\
&  =  A_0(k,\ell)
           \left(
                      \frac{4}{B}
                         \log \left(\frac{Be^\gamma k}{4}\right) - \nu + o(1) \right).
\end{align*}

This is positive if
\begin{equation}
\frac{B e^\gamma k }{4} \geq e^{\frac{B\nu }{4}}(1 + o(1)),
\label{eq:3.11}
\end{equation}
and this  proves Theorem~\ref{th:1}.

We remark that in the above proof, we are finding ``unsifted'' $E_2$-numbers; 
i.e., the $E_2$-numbers found in the proof can have small prime factors. 
However, it should be clear from the argument, that if one desires, one may take 
$Y$ in the definition of $\beta$ to be any function of $N$ such that $\log Y/\log N\to 0$
as $N\to \infty$, and the same argument goes through.

In order to show Corollary~\ref{cor:1}, we have only to note
that if $p_1 < p_2 < \ldots $ are the consecutive primes then
\begin{equation}
\mathcal H = \{ p_{\pi(k) + 1} \dots p_{\pi(k) + k} \}
\label{eq:3.12}
\end{equation}
forms an admissible $k$-tuple and $p_{\pi(k)+k} \sim k \log k$.

Now we consider Theorem \ref{th:2}. Let $\cS$ be as defined in \eqref{eq:cSDef} with $\nu=1$.
By Theorems \ref{T:S0Est} and \ref{T:S1Est}, we see that 
\begin{equation*}
\cS \asymptotic \frac{\fS(\cL) N(\log R)^{k} } {(k-1)!} J,
\end{equation*}
where 
\begin{equation*}
J= \frac{k(k-1)}{B} (J_1+J_2+J_3) - J_0.
\end{equation*}
We take $k=3$, $B=4$, $\eta=1/144$, and $P(x)=1+6x$. 
Straightforward computations show that
\begin{align*}
     J_0= & \frac{38}{15} = 2.5333,\ldots ,\\
     J_1= & 4824\log\left(\frac{143}{108}\right) - 
                 \frac{13641020155}{10077696} 
              = 0.57625 \ldots,\\
     J_2= & -\frac{77824}{15}\log\left(\frac{143}{108}\right) + 
     \frac{14680965985}{10077696} =
	       0.36202 \ldots, \\
     J_3= & \frac{41}{60}\log 3 = 0.75071 \ldots, \\
   J = & 
    \frac{41}{40}\log 3 
     -\frac{2732}{5}\log\left(\frac{143}{108}\right) + 
        \frac{852438101}{5598720} 
         = 
     0.00016493 \ldots .
\end{align*}

For Theorem 3, we take
 $k=2, B=4,\eta=1/10$, and consider the sum
\begin{equation*}
\ctS(\cL) = 
\sum_{N < n \le 2N} 
 \left\{ \sum_{j=1}^2  ( \beta(L_j(n)) + \varpi(L_j(n)) - 1\right\}
     \left(\sum_{d|P_{\cL}(n)} \lambda_d\right)^2.
\end{equation*}
From Theorems \ref{T:S0Est}, \ref{T:S1Est}, and \ref{th:9}, we see that
\begin{equation*}
\ctS(\cL) \asymptotic N\fS(\cL)(\log R)^2 J,
\end{equation*}
where 
\begin{equation*}
J= \frac{1}{2} (J_1+J_2+J_3+J_{\varpi}) - J_0.
\end{equation*}
With $P(x)=1+x$, we see that 
\begin{align*}
     J_0= & \frac{11}{12} =  0.91667  \ldots ,\\
     J_1= & -144 \log(6/5) + \frac{66363}{2500}= 0.29089 \ldots, \\
     J_2= & \frac{2048}{15} \log(6/5) - \frac{308429}{12500} = 0.21864 \ldots, \\
     J_3= & \frac{19}{30}\log 3 = 0.69578 \ldots, \\
     J_{\varpi}= & \frac{19}{30} = 0.63333 \ldots,\\
     J    = & \frac{19}{60} \log 3 - \frac{56}{15} \log (6/5) + \frac{4193}{12500} = 0.00266 \ldots,
\end{align*}
and the theorem follows. The result of \eqref{eq:1.21} follows by taking $\cL= \{n, n-d\}$.

Now we mention the slight changes which lead to the proofs of
Theorems~\ref{th:4}--\ref{th:6}.

Theorem~\ref{th:4} follows from the proof of Theorem~\ref{th:1}
by taking $B = 60$ in view of \eqref{eq:1.28}--\eqref{eq:1.29}.
For Theorem~\ref{th:5} we have to restrict $p$ and $q$ to primes
of the form $4m + 1$.
This means that the density of both $p$ and $q$ is half of that
of all primes, therefore we obtain finally for all $S_i$ and
$J_i$ $(i = 2, 4, 5, 6)$ a quantity which is $1/4$ of that in
the proof of Theorem~\ref{th:1}, which has the same effect as to
writing $4\nu $ in place of~$\nu $.

Finally, the proof of Theorem~\ref{th:6} is just a combination
of the proofs of Theorems~\ref{th:4} and \ref{th:5}.
The result is that we have to take $B = 60$ as in
Theorem~\ref{th:4} and to replace $\nu $ by $4 \nu $ 
as in Theorem~\ref{th:5}.
This leads finally to \eqref{eq:1.37}.

Corollaries \ref{cor:3}--\ref{cor:5} follow from
Theorems~\ref{th:4}--\ref{th:6} in the same way as
Corollary~\ref{cor:1} follows from Theorem~\ref{th:1} 
(see \eqref{eq:3.12}).

\bigskip
\footnotesize
D. A. Goldston,
Department of Mathematics,
San Jose
State University,
San Jose, CA 95192,
USA,
e-mail: goldston@math.sjsu.edu

S. W. Graham,
Department of Mathematics,
Central Michigan University,
Mt. Pleasant, MI 48859, USA,
email: graha1sw@cmich.edu

J. Pintz,
R\'enyi Mathematical Institute of the Hungarian Academy
of Sciences,
H-1364 Budapest,
P.O.B. 127,
Hungary,
e-mail: pintz@renyi.hu

C. Y. Y{\i}ld{\i}r{\i}m,
Department of Mathematics,
Bo\~{g}azi\c{c}i University,
Istanbul 34342,
\& \newline 
Feza G\"{u}rsey Enstit\"{u}s\"{u}, \c{C}engelk\"{o}y, Istanbul,
P.K. 6, 81220, Turkey,
e-mail: yalciny@boun.edu.tr

\end{document}